\documentclass[letterpaper]{article}
\usepackage{aaai}
\usepackage{times}
\usepackage{helvet}
\usepackage{courier}
\usepackage{amssymb}
\usepackage{amsmath}
\usepackage{graphicx}
\usepackage{amsthm}
\usepackage[framed,numbered]{mcode}
\usepackage{multirow}
\usepackage{bm}
\usepackage{algorithm}
\usepackage{algorithmic}
\usepackage{enumerate}
\usepackage{subcaption}

\newtheorem{theo}{Theorem}

\begin{document}
%

\title{Proximal Iteratively Reweighted Algorithm with Multiple Splitting for Nonconvex Sparsity Optimization}

\author{Canyi Lu$^1$, Yunchao Wei$^2$, Zhouchen Lin$^{3,}$\thanks{Corresponding author.}, Shuicheng Yan$^1$\\
$^1$ Department of Electrical and Computer Engineering, National University of Singapore\\
$^2$ Institute of Information Science, Beijing Jiaotong University\\
$^3$ Key Laboratory of Machine Perception (MOE), School of EECS, Peking University\\
{\tt\small canyilu@gmail.com, wychao1987@gmail.com, zlin@pku.edu.cn, eleyans@nus.edu.sg}
}

\maketitle
\begin{abstract}
\begin{quote}
This paper proposes the Proximal Iteratively REweighted (PIRE) algorithm for solving a general problem, which involves a large body of nonconvex sparse and structured sparse related problems. Comparing with previous iterative solvers for nonconvex sparse problem, PIRE is much more general and efficient. The computational cost of PIRE in each iteration is usually as low as the state-of-the-art convex solvers. We further propose the PIRE algorithm with Parallel Splitting (PIRE-PS) and PIRE algorithm with Alternative Updating (PIRE-AU) to handle the multi-variable problems. In theory, we prove that our proposed methods converge and any limit solution is a stationary point. Extensive experiments on both synthesis and real data sets demonstrate that our methods achieve comparative learning performance, but are much more efficient, by comparing with previous nonconvex solvers.
\end{quote}
\end{abstract}

\section{Introduction}
\label{sec:intro}
This paper aims to solve the following general problem
\begin{equation}\label{eq_genp}
\min_{\mathbf{x}\in\mathbb{R}^n} F(\mathbf{x})=\lambda f(\mathbf{g}(\mathbf{x}))+h(\mathbf{x}),
\end{equation}
where $\lambda>0$ is a parameter, and the functions in the above formulation satisfy the following conditions:
\begin{itemize}
\item[\textbf{C1}] $f(\mathbf{y})$ is nonnegative, concave and increasing.
\item[\textbf{C2}] $\mathbf{g}(\mathbf{x}): \mathbb{R}^n\rightarrow{\mathbb{R}^d}$ is a nonnegative multi-dimensional function, such that the following problem
\begin{equation}\label{eq_proxw}
\min_{\mathbf{x}\in\mathbb{R}^n}\lambda\langle\mathbf{w},\mathbf{g}(\mathbf{x})\rangle+\frac{1}{2}||\mathbf{x}-\mathbf{b}||_2^2,
\end{equation}
is convex and can be cheaply solved for any given nonnegative $\mathbf{w}\in\mathbb{R}^d$.
\item[\textbf{C3}] $h(\mathbf{x})$ is a smooth function of type $C^{1,1}$, i.e., continuously differentiable with the Lipschitz continuous gradient
\begin{equation}
||\nabla h(\mathbf{x})-\nabla h(\mathbf{y})||\leq L(h)||\mathbf{x}-\mathbf{y}|| \ \text{ for any } \ \mathbf{x}, \mathbf{y}\in\mathbb{R}^{ n},
\end{equation}
$L(h)>0$ is called the Lipschitz constant of $\nabla h$.
\item[\textbf{C4}] $\lambda f(\mathbf{g}(\mathbf{x}))+h(\mathbf{x})\rightarrow\infty$ iff $||\mathbf{x}||_2\rightarrow\infty$.
\end{itemize}

Note that problem (\ref{eq_genp}) can be convex or nonconvex. Though $f(\mathbf{y})$ is concave, $f(\mathbf{g}(\mathbf{x}))$ can be convex w.r.t $\mathbf{x}$. Also $f(\mathbf{y})$ and $\mathbf{g}(\mathbf{x})$ are not necessarily smooth, and $h(\mathbf{x})$ is not necessarily convex.

Based on different choices of $f$, $\mathbf{g}$, and $h$, the general problem (\ref{eq_genp}) involves many sparse representation models, which have many important applications in machine learning and computer vision \cite{SRC,beck2009fast,jacob2009group,gong2012robust}. For the choice of $h$, the least square and logistic loss functions are two most widely used ones which satisfy (\textbf{C3}):
\begin{equation}
h(\mathbf{x})=\frac{1}{2}||\mathbf{A}\mathbf{x}-\mathbf{b}||_2^2, \text{or} \ \frac{1}{n}\sum_{i=1}^{n}\log(1+\exp(-b_i\textbf{a}_i^T\textbf{x})),
\end{equation}
where $\textbf{A}=[\mathbf{a}_1^T;\cdots;\mathbf{a}_n^T]\in\mathbb{R}^{n\times d}$, and $\mathbf{b}\in\mathbb{R}^n$. As for the choice of $\mathbf{g}(\mathbf{x})$, $|\mathbf{x}|$ (absolute value of $\mathbf{x}$ element-wise) and $\mathbf{x}^2$ (square of $\mathbf{x}$ element-wise) are widely used. One may also use $\mathbf{g}(\mathbf{X})=||\mathbf{x}_i||_2$ ($\mathbf{x}_i$ denotes the $i$-th column of $\mathbf{X}$) when pursuing column sparsity of a matrix $\mathbf{X}$. As for the choice of $f$, almost all the existing nonconvex surrogate functions of the $\ell_0$-norm are concave on $(0,\infty)$. In element-wise, they include $\ell_p$-norm $y^p$ ($0<p<1$) \cite{knight2000asymptotics}, logarithm function $\log(y)$ \cite{candes2008enhancing}, smoothly clipped absolute deviation \cite{fan2001variable}, and minimax concave penalty \cite{zhang2010nearly}. 

The above nonconvex penalties can be further extended to define structured sparsity \cite{jacob2009group}. For example, let $\mathbf{x}=[\mathbf{x}_1;\cdots;\mathbf{x}_G]$. By taking $g(\mathbf{x})=[||\mathbf{x}_1||_2;\cdots;||\mathbf{x}_G||_2]$ and $f(\mathbf{y})=\sum_if_i(y_i)$, with $f_i$ being any of the above concave functions, then $f(\mathbf{g}(\mathbf{x}))$ is the nonconvex group Lasso $\sum_{i}f_i(||\mathbf{x}_i||_2)$. By taking $f(\mathbf{y})=\sum_iy_i$, $f(\mathbf{g}(\mathbf{x}))=\sum_i||\mathbf{x}_i||_2$ is the group Lasso. 

Problem (\ref{eq_genp}) contains only one variable. We will show that our proposed model can be naturally used for handling problem with several variables (which we mean a group of variables that can be updated simultaneously due to the separability structure of the problem). An example for multi-task learning can be found in \cite{gong2012robust}.

\subsection{Related Works}

If the condition (\textbf{C3}) holds and
\begin{equation}\label{eq_prom}
\min_{\mathbf{x}} \lambda f(\mathbf{g}(\mathbf{x}))+\frac{1}{2}||\mathbf{x}-\mathbf{b}||_2^2,
\end{equation}
can be cheaply computed, then problem (\ref{eq_genp}) can be solved by iteratively solving a series of problem (\ref{eq_prom}) \cite{gong2013general}. Such an updating procedure is the same as the ISTA algorithm \cite{beck2009fast}, which is originally for convex optimization. It can be proved that any accumulation point of $\{\mathbf{x}^k\}$ is a stationary point of problem (\ref{eq_genp}). If $f(\mathbf{g}(\mathbf{x}))$ and $h(\mathbf{x})$ are convex, the Fast ISTA algorithm (FISTA) \cite{beck2009fast} converges to the globally optimal solution with a convergence rate $O(1/T^2)$ ($T$ is the iteration number). But for nonconvex optimization, it is usually very difficult to get the globally optimal solution to problem (\ref{eq_prom}). Sometimes, it is also not easy even if $f(\mathbf{g}(\mathbf{x}))$ is convex.

The multi-stage algorithm in \cite{zhang2008multi} solves problem (\ref{eq_genp}) by solving a series of convex problem.
\begin{equation}\label{eq_mutipro}
\min_{\mathbf{x}}\lambda\langle\mathbf{w}^k,\mathbf{g}(\mathbf{x})\rangle+h(\mathbf{x}).
\end{equation}
However, solving such a convex problem requires other iterative solvers which is not efficient. It also fails when $h(\mathbf{x})$ is nonconvex.

More specially, the Iteratively Reweighted L1 (IRL1) \cite{chenconvergence} and Iteratively Reweighted Least Squares (IRLS) \cite{lai2013improved} algorithms are special cases of the multi-stage algorithm. They aim to solve the following $\ell_p$-regularization problem
\begin{equation}\label{eq_lp}
\min_{\mathbf{x}}  \lambda||\mathbf{x}||_p^p+\frac{1}{2}||\mathbf{A}\mathbf{x}-\mathbf{b}||_2^2.
\end{equation}
The above problem is NP-hard. IRL1 instead considers the following relaxed problem
\begin{equation}\label{eq_smo11}
\min_{\mathbf{x}} \lambda\sum_{i=1}^n(|x_i|+\epsilon)^p+\frac{1}{2}||\mathbf{A}\mathbf{x}-\mathbf{b}||_2^2,
\end{equation}
with $0<\epsilon\ll1$. IRL1 updates $\mathbf{x}^{k+1}$ by solving
\begin{equation}\label{eq_IRL1w}
\mathbf{x}^{k+1}=\arg\min_{\mathbf{x}}\lambda \sum_{i=1}^nw^k_i|x_i|+\frac{1}{2}||\mathbf{A}\mathbf{x}-\mathbf{b}||_2^2,
\end{equation}
with $w_i^k={p}/{(|x_i^k|+\epsilon)^{1-p}}$. Problem (\ref{eq_smo11}) is a special case of (\ref{eq_genp}) by letting $f(\mathbf{y})=\sum_i(y_i+\epsilon)^p$ ($0<p<1$) and $\mathbf{g}(\mathbf{x})=|\mathbf{x}|$. However, IRL1 is not efficient since it has to solve a number of nonsmooth problem (\ref{eq_IRL1w}) by using some other convex optimization methods, e.g. FISTA.

The other method, IRLS, smooths problem (\ref{eq_lp}) as
\begin{equation}\label{eq_IRLSp}
\min \lambda\sum_{i=1}^n(x_i^2+\epsilon)^{\frac{p}{2}}+\frac{1}{2}||\mathbf{A}\mathbf{x}-\mathbf{b}||_2^2,
\end{equation}
and updates $\mathbf{x}^{k+1}$ by solving
\begin{equation}\label{eq_sol2}
{\lambda}\text{Diag}({\mathbf{w}^k})\mathbf{x}+\mathbf{A}^T(\mathbf{A}\mathbf{x}-\mathbf{b})=0,
\end{equation}
with $w_i^k={p}/{\left((x_i^k)^2+\epsilon\right)^{1-\frac{p}{2}}}$. Problem (\ref{eq_IRLSp}) is also a special case of (\ref{eq_genp}) by taking $f(\mathbf{y})=\sum_{i}(y_i+\epsilon)^{\frac{p}{2}}$ and $\mathbf{g}(\mathbf{x})=\mathbf{x}^2$. However, the obtained solution by IRLS may not be naturally sparse, or it may require a lot of iterations to get a sparse solution. One may perform thresholding appropriately to achieve a sparse solution, but there is no theoretically sound rule to suggest a correct threshold.

Another related work is \cite{lu2012iterative} which aims to solve
\begin{equation}\label{eq_smo1sss}
\min_{\mathbf{x}} \lambda\sum_{i=1}^n(|x_i|+\epsilon)^p+h(\mathbf{x}).
\end{equation}
In each iteration, $\mathbf{x}$ is efficiently updated by solving a series of problem
\begin{equation}\label{eq_wl1prx}
\min_{\mathbf{x}}\lambda\langle\mathbf{w}^k,|\mathbf{x}|\rangle+\frac{1}{2}||\mathbf{x}-\mathbf{b}||_2^2.
\end{equation}
But their solver is only for problem (\ref{eq_smo1sss}) which is a special case of (\ref{eq_genp}). The convergence proofs also depend on the special property of the $\ell_p$-norm, thus is not general.

Furthermore, previous iterative algorithms can only solve the problem with only one variable. They cannot be naively generalized to solve multi-variable problems. However, there are many problems involving two or more variables, e.g. stable robust principle component analysis \cite{zhou2010stable} and robust multi-task feature learning \cite{gong2012robust}. So it is desirable to extend the iteratively reweighted algorithms for the multi-variable case.

%
\subsection{Contributions}
In this work, we propose a novel method to solve the general problem (\ref{eq_genp}), and address the scalablity and multi-variable issues. In each iteration we only need to solve problem (\ref{eq_proxw}), whose computational cost is usually the same as previous state-of-the-art first-order convex methods. This method is named as Proximal Iteratively REweighted (PIRE) algorithm. We further propose two multiple splitting versions of PIRE: PIRE with Parallel Splitting (PIRE-PS) and PIRE with Alternative Updating (PIRE-AU) to handle the multi-variable problem. Parallel splitting makes the algorithm highly parallelizable, making PIRE-PS suitable for distributed computing. This is important for large scale applications. PIRE-AU may converge faster than PIRE-PS. In theory, we prove that any sequences generated by PIRE, PIRE-PS and PIRE-AU are bounded and any accumulation point is a stationary point. To the best of our knowledge, PIRE-PS and PIRE-AU are the first two algorithms for problem (\ref{eq_genp}) with multi-variables. If problem (\ref{eq_genp}) is convex, the obtained solution is globally optimal.




\section{Proximal Iteratively Reweighted Algorithm}
In this section, we show how to solve problem (\ref{eq_genp}) by our Proximal Iteratively Reweighted (PIRE) algorithm. Instead of minimizing $F(\mathbf{x})$ in (\ref{eq_genp})  directly, we update $\mathbf{x}^{k+1}$ by minimizing the sum of two surrogate functions, which correspond to two terms of $F(\mathbf{x})$, respectively.

\begin{algorithm}[t]
\caption{Solving problem (\ref{eq_genp}) by PIRE}
\textbf{Input:} $\mu>\frac{L(h)}{2}$, where $L(h)$ is the Lipschitz constant of $h(\mathbf{x})$.  \\
\textbf{Initialize:} $k=0$, $\mathbf{w}^k$.\\
\textbf{Output:} $\mathbf{x}^*$. \\
\textbf{while} not converge \textbf{do}
\begin{enumerate}
  \item Update $\mathbf{x}^{k+1}$ by solving the following problem \\
  \begin{equation*}
  \begin{split}
  \mathbf{x}^{k+1}=&\arg\min_{\mathbf{x}}{\lambda}\langle\mathbf{w}^k,\mathbf{g}(\mathbf{x})\rangle+\frac{\mu}{2}\left\|\mathbf{x}-\left(\mathbf{x}^{k}-\frac{1}{\mu}\nabla h(\mathbf{x}^k) \right)\right\|^2.
  \end{split}
  \end{equation*}
  \item Update the weight $\mathbf{w}^{k+1}$ by  \\
    \begin{equation*}
\mathbf{w}^{k+1}\in-\partial\left(-f(\mathbf{g}(\mathbf{x}^{k+1}))\right).
\end{equation*}
\textbf{end while}
\end{enumerate}
\label{alg_PIRE}
\end{algorithm}

First, note that $f(\mathbf{y})$ is concave, $-f(\mathbf{y})$ is convex. By the definition of subgradient of the convex function, we have
\begin{equation}\label{eq_lpmaj0}
-f(\mathbf{g}(\mathbf{x}))\geq -f(\mathbf{g}(\mathbf{x}^k))+\langle -\mathbf{w}^k,\mathbf{g}(\mathbf{x})-\mathbf{g}(\mathbf{x}^k)\rangle,
\end{equation}
where $-\mathbf{w}^k$ is the subgradient of $-f(\mathbf{y})$ at $\mathbf{y}=\mathbf{g}(\mathbf{x}^k)$, i.e.
\begin{equation}\label{updatew}
-\mathbf{w}^k\in\partial\left(-f(\mathbf{g}(\mathbf{x}^k))\right) \ \text{or} \  \mathbf{w}^k\in-\partial\left(-f(\mathbf{g}(\mathbf{x}^k))\right).
\end{equation}
Eqn (\ref{eq_lpmaj0}) is equivalent to
\begin{equation}\label{eq_lpmaj}
f(\mathbf{g}(\mathbf{x}))\leq f(\mathbf{g}(\mathbf{x}^k))+\langle \mathbf{w}^k,\mathbf{g}(\mathbf{x})-\mathbf{g}(\mathbf{x}^k)\rangle.
\end{equation}
Then $f(\mathbf{g}(\mathbf{x}^k))+\langle \mathbf{w}^k,\mathbf{g}(\mathbf{x})-\mathbf{g}(\mathbf{x}^k)\rangle$ is used as a surrogate function of $f(\mathbf{g}(\mathbf{x}))$.

The loss function $h(\mathbf{x})$, which has Lipschitz continuous gradient, owns the following property \cite{bertsekas1995nonlinear}
\begin{equation}\label{eq_lipp}
h(\mathbf{x})\leq h(\mathbf{y})+\langle \nabla h(\mathbf{y}),\mathbf{x}-\mathbf{y}\rangle+\frac{L(h)}{2}||\mathbf{x}-\mathbf{y}||_2^2.
\end{equation}
Let $\mathbf{y}=\mathbf{x}^k$, $h(\mathbf{x}^k)+\langle \nabla h(\mathbf{x}^k),\mathbf{x}-\mathbf{x}^k\rangle+\frac{L(h)}{2}||\mathbf{x}-\mathbf{x}^k||_2^2$ is used as a surrogate function of $h(\mathbf{x})$.

Combining (\ref{eq_lpmaj}) and (\ref{eq_lipp}), we update $\mathbf{x}^{k+1}$ by minimizing the sum of these two surrogate functions
\begin{equation}\label{eq_updatex}
\begin{split}
&\mathbf{x}^{k+1}\\
=&\arg\min_{\mathbf{x}} f(\mathbf{g}(\mathbf{x}^k))+\langle \mathbf{w}^k,\mathbf{g}(\mathbf{x})-\mathbf{g}(\mathbf{x}^k)\rangle\\
&+h(\mathbf{x}^k)+\langle \nabla h(\mathbf{x}^k),\mathbf{x}-\mathbf{x}^k\rangle+\frac{\mu}{2}||\mathbf{x}-\mathbf{x}^k||_2^2\\
=&\arg\min_{\mathbf{x}} \lambda\langle \mathbf{w}^k,\mathbf{g}(\mathbf{x})\rangle+\frac{\mu}{2}\left\|\mathbf{x}-\left(\mathbf{x}^k-\frac{1}{\mu}\nabla h(\mathbf{x}^k)\right)\right\|_2^2,
\end{split}
\end{equation}
where $\mathbf{w}^k$ is also called the weight corresponding to $\mathbf{g}(\mathbf{x}^k)$.

For the choice of $\mu$ in (\ref{eq_updatex}), our theoretical analysis shows that $\mu>L(h)/2$ guarantees the convergence of the proposed algorithm. Note that $f$ is concave and increasing, this guarantees that $\mathbf{w}^k$ in (\ref{updatew}) is nonnegative. Usually problem (\ref{eq_updatex}) can be cheaply computed based on the condition (\textbf{C2}). For example, if $\mathbf{g}(\mathbf{x})=|\mathbf{x}|$, solving problem (\ref{eq_updatex}) costs only $O(n)$.
%
Such computational cost is the same as the state-of-the-art convex solvers for $\ell_1$-minimization. This idea leads to the Proximal Iteratively REweighted (PIRE) algorithm, as shown in Algorithm \ref{alg_PIRE}. In the next section, we will prove that the sequence generated by PIRE is bounded and any accumulation point is a stationary point of problem (\ref{eq_genp}).

\section{Convergence Analysis of PIRE}

\begin{theo}\label{thm_pro}
Let $D=F(\mathbf{x}^1)$, and $\mu>\frac{L(h)}{2}$, where $L(h)$ is the Lipschitz constant of $h(\mathbf{x})$. The sequence $\{\mathbf{x}^k\}$ generated in Algorithm \ref{alg_PIRE} satisfies the following properties:
\begin{enumerate}[(1)]
\item $F(\mathbf{x}^k)$ is monotonically decreasing. Indeed,
\begin{equation*}
F(\mathbf{x}^k)-F(\mathbf{x}^{k+1})\geq\left(\mu-\frac{L(h)}{2}\right)||\mathbf{x}^k-\mathbf{x}^{k+1}||^2;
\end{equation*}
\item The sequence $\{\mathbf{x}^k\}$ is bounded;
\item $\sum\limits_{k=1}^\infty||\mathbf{x}^k-\mathbf{x}^{k+1}||_F^2\leq \frac{2D}{2\mu-L(h)}$. In particular, we have $\lim\limits_{k\rightarrow\infty}(\mathbf{x}^k-\mathbf{x}^{k+1})=\bm{0}$.
\end{enumerate}
\end{theo}
\textbf{Proof}. Since $\mathbf{x}^{k+1}$ is the globally optimal solution to problem (\ref{eq_updatex}), the zero vector is contained in the subgradient with respect to $\mathbf{x}$. That is, there exists $\mathbf{v}^{k+1}\in\partial\langle\mathbf{w}^k,\mathbf{g}(\mathbf{x}^{k+1})\rangle$ such that
\begin{equation}\label{eq_proof1}
\lambda\mathbf{v}^{k+1}+\nabla h(\mathbf{x}^k)+\mu(\mathbf{x}^{k+1}-\mathbf{x}^{k})=\mathbf{0}.
\end{equation}
A dot-product with $\mathbf{x}^{k+1}-\mathbf{x}^k$ on both sides of (\ref{eq_proof1}) gives
\begin{equation}\label{eq_proof2}
\begin{split}
&\lambda\left\langle \mathbf{v}^{k+1},\mathbf{x}^{k+1}-\mathbf{x}^k\right\rangle+\left\langle \nabla h(\mathbf{x}^k),\mathbf{x}^{k+1}-\mathbf{x}^k\right\rangle\\
&+\mu||\mathbf{x}^{k+1}-\mathbf{x}^k||^2=0.
\end{split}
\end{equation}
Recalling the definition of the subgradient of the convex function, we have
\begin{equation}\label{eq_proof3}
\begin{split}
\langle\mathbf{w}^k,\mathbf{g}(\mathbf{x}^k)-\mathbf{g}(\mathbf{x}^{k+1})\rangle\geq\left\langle\mathbf{v}^{k+1},\mathbf{x}^k-\mathbf{x}^{k+1}\right\rangle.
\end{split}
\end{equation}
Combining (\ref{eq_proof2}) and (\ref{eq_proof3}) gives
\begin{equation}\label{eq_proof4}
\begin{split}
&\lambda\langle\mathbf{w}^k,\mathbf{g}(\mathbf{x}^k)-\mathbf{g}(\mathbf{x}^{k+1})\rangle\\
\geq& -\left\langle \nabla h(\mathbf{x}^k),\mathbf{x}^k-\mathbf{x}^{k+1}\right\rangle+\mu||\mathbf{x}^{k+1}-\mathbf{x}^k||^2.
\end{split}
\end{equation}
Since $f$ is concave, similar to $(\ref{eq_lpmaj})$, we get
\begin{equation}\label{eq_proof5}
f(\mathbf{g}(\mathbf{x}^k))-f(\mathbf{g}(\mathbf{x}^{k+1}))\geq\langle\mathbf{w}^k,\mathbf{g}(\mathbf{x}^k)-\mathbf{g}(\mathbf{x}^{k+1})\rangle.
\end{equation}
By the condition (\textbf{C3}), we have
\begin{equation}\label{eq_proof6}
\begin{split}
&h(\mathbf{x}^k)-h(\mathbf{x}^{k+1})\\
\geq & \left\langle \nabla h(\mathbf{x}^k),\mathbf{x}^k-\mathbf{x}^{k+1}\right\rangle-\frac{L(h)}{2}||\mathbf{x}^{k+1}-\mathbf{x}^k||^2.
\end{split}
\end{equation}
Now, combining (\ref{eq_proof4})(\ref{eq_proof5}) and (\ref{eq_proof6}), we have
\begin{equation}\label{eq_proof66}
\begin{split}
&F(\mathbf{x}^k)-F(\mathbf{x}^{k+1})\\
=&\lambda f(\mathbf{g}(\mathbf{x}^k))-\lambda f(\mathbf{g}(\mathbf{x}^{k+1}))+h(\mathbf{x}^k)-h(\mathbf{x}^{k+1})\\
\geq&\left(\mu-\frac{L(h)}{2}\right)||\mathbf{x}^{k+1}-\mathbf{x}^k||^2\geq0.
\end{split}
\end{equation}
Hence $F(\mathbf{x}^k)$ is monotonically decreasing. Summing all the above inequalities for $k\geq 1$, it follows that
\begin{equation}\label{eq_proof7}
D=F(\mathbf{x}^1)\geq\left(\mu-\frac{L(h)}{2}\right)\sum_{k=1}^{\infty}||\mathbf{x}^{k+1}-\mathbf{x}^k||^2.
\end{equation}
This implies that $\lim\limits_{k\rightarrow\infty}(\mathbf{x}^{k+1}-\mathbf{x}^k)=\mathbf{0}$. Also $\{\mathbf{x}^k\}$ is bounded due to the condition (\textbf{C4}).
$\hfill\blacksquare$

\begin{theo}\label{thm_con}
Let $\{\mathbf{x}^k\}$ be the sequence generated in Algorithm \ref{alg_PIRE}. Then any accumulation point of $\{\mathbf{x}^k\}$ is a stationary point $\mathbf{x}^*$ of problem (\ref{eq_genp}). Furthermore, for every $n\geq1$, we have
\begin{equation}\label{eq_converrate}
\min_{1\leq k\leq n}||\mathbf{x}^{k+1}-\mathbf{x}^k||_2^2\leq\frac{F(\mathbf{x}^1)-F(\mathbf{x}^*)}{n\left(\mu-\frac{L(h)}{2}\right)}.
\end{equation}
\end{theo}

Please refer to the Supplementary Material for the proof.

We conclude this section with the following remarks:
\begin{enumerate}[(1)]
\item When proving the convergence of IRL1 for solving problem (\ref{eq_smo11}) or (\ref{eq_smo1sss}) in \cite{chenconvergence,lu2012iterative}, they use the Young's inequality which is a special property of the function $y^p$ ($0<p<1$)
\begin{equation}\label{eq_proof55555}
\sum_{i=1}^n(|x_i^k|+\epsilon)^p-(|x_i^{k+1}|+\epsilon)^p\geq\sum_{i=1}^{n}w^k_i\left(\left|x^{k}_i\right|-\left|x^{k+1}_i\right|\right),
\end{equation}
where $w^k_i=p/(|x_i^k|+\epsilon)^{1-p}$. Eqn (\ref{eq_proof55555}) is a special case of (\ref{eq_proof5}). But (\ref{eq_proof5}) is obtained by using the concavity of $f(\mathbf{y})$, which is much more general.
\item In Eqn (\ref{eq_converrate}), $||\mathbf{x}^{k+1}-\mathbf{x}^k||_2$ is used to measure the convergence rate of the algorithm. The reason is that $||\mathbf{x}^{k+1}-\mathbf{x}^k||_2\rightarrow{0}$ is a necessary optimality condition as shown in the Theorem \ref{thm_pro}.
\item PIRE requires that $\mu>L(h)/2$. But sometimes the Lipschitz constant $L(h)$ is not known, or it is not computable for large scale problems. One may use the backtracking rule to estimate $\mu$ in each iteration \cite{beck2009fast}. PIRE with multiple splitting shown in the next section also eases this problem.
\end{enumerate}

\section{PIRE with Multiple Splitting}
In this section, we will show that PIRE can also solve multi-variable problem as follows
\begin{equation}\label{eq_mspro}
\min_{\mathbf{x}_1,\cdots,\mathbf{x}_S} F(\mathbf{x})=\lambda\sum_{s=1}^Sf_s(\mathbf{g}_s(\mathbf{x}_s))+h(\mathbf{x}_1,\cdots,\mathbf{x}_S),
\end{equation}
where $f_s$ and $\mathbf{g}_s$ holds the same assumptions as $f$ and $\mathbf{g}$ in problem (\ref{eq_genp}), respectively. Problem (\ref{eq_mspro}) is similar to problem (\ref{eq_genp}), but splits $\mathbf{x}$ into $\mathbf{x}=[\mathbf{x}_1;\cdots ;\mathbf{x}_S]\in\mathbb{R}^n$, where $\mathbf{x}_s\in\mathbb{R}^{n_s}$, and $\sum_{i=1}^Sn_s=n$.


Based on different assumptions of $h(\mathbf{x}_1,\cdots,\mathbf{x}_S)$, we have two splitting versions of the PIRE algorithm. They use different updating orders of the variables.


\subsection{PIRE with Parallel Splitting}
If we still assume that (\textbf{C3}) holds, i.e. $h(\mathbf{x}_1,\cdots,\mathbf{x}_S)$ has a Lipschitz continuous gradient, with Lipschitz constant $L(h)$, PIRE is naturally parallelizable. In each iteration, we parallelly update $\mathbf{x}_s^{k+1}$ by
  \begin{equation}\label{updatexinps}
  \begin{split}
  \mathbf{x}_s^{k+1}=&\arg\min_{\mathbf{x_s}}{\lambda}\langle\mathbf{w}_s^k,\mathbf{g}_s(\mathbf{x}_s)\rangle\\
  &+\frac{\mu}{2}\left\|\mathbf{x}_s-\left(\mathbf{x}^{k}_s-\frac{1}{\mu}\nabla_s h\left(\mathbf{x}^k_1,\cdots,\mathbf{x}^k_S\right) \right)\right\|^2,
  \end{split}
  \end{equation}
where the notion $\nabla_s h\left(\mathbf{x}_1,\cdots,\mathbf{x}_S\right)$ denotes the gradient w.r.t $\mathbf{x}_s$, $\mu>L(h)/2$, and $\mathbf{w}_s^k$ is the weight vector corresponding to $\mathbf{g}(\mathbf{x}_s^{k})$, which can be computed by
\begin{equation}\label{eq_updatewinpa}
\mathbf{w}_s^k\in-\partial\left(-f_s(\mathbf{g}_s(\mathbf{x}_s^k))\right), \ s=1,\cdots,S.
\end{equation}
When updating $\mathbf{x}_s$ in the $(k+1)$-th iteration, only the variables in the $k$-th iteration are used. Thus the variables $\mathbf{x}_s^{k+1}$, $s=1,\cdots,S$, can be updated in parallel. This is known as Jacobi iteration in numerical algebra \cite{LADMPSAP}. This algorithm is named as PIRE with Parallel Splitting (PIRE-PS). Actually the updating rule of PIRE-PS is the same as PIRE, but in parallel. It is easy to check that the proofs in Theorem \ref{thm_pro} and \ref{thm_con} also hold for PIRE-PS.

For some special cases of $h(\mathbf{x}_1,\cdots,\mathbf{x}_S)$, we can use different $\mu_s$, usually smaller than $\mu$, for updating $\mathbf{x}_s^{k+1}$. This may lead to faster convergence \cite{zuo2011generalized}. If $h(\mathbf{x}_1,\cdots,\mathbf{x}_S)=\frac{1}{2}\left\|\sum_{s=1}^S\mathbf{A}_s\mathbf{x}_s-\mathbf{b}\right\|_2^2$, we can update $\mathbf{x}_s^{k+1}$ by
\begin{equation}\label{eq_psupx}
  \begin{split}
  \mathbf{x}_s^{k+1}=&\arg\min_{\mathbf{x_s}}\lambda\langle\mathbf{w}_s^k,\mathbf{g}_s(\mathbf{x}_s)\rangle\\
  &+\frac{\mu_s}{2}\left\|\mathbf{x}_s-\left(\mathbf{x}^{k}_s-\frac{1}{\mu_s}\mathbf{A}_s^T(\mathbf{A}\mathbf{x}^k-\mathbf{b}\right)\right\|_2^2,
  \end{split}
\end{equation}
where $\mu_s>L_s(h)/2$ and $L_s(h)=||\mathbf{A}_s||_2^2$ is the Lipschitz constant of $\nabla_sh(\mathbf{x}_1,\cdots,\mathbf{x}_S)$. If the size of $\mathbf{A}$ is very large, $L(h)=||\mathbf{A}||_2^2$ may not be computable. We can split it to $\mathbf{A}=[\mathbf{A}_1,\cdots,\mathbf{A}_S]$, and compute each $L_s(h)=||\mathbf{A}_s||_2^2$ instead. Similar convergence results in Theorem \ref{thm_pro} and \ref{thm_con} also hold by updating $\mathbf{x}_s^{k+1}$ in (\ref{eq_psupx}). For detailed proofs, please refer to the Supplementary Material. A main difference of the convergence poof is that we use the Pythagoras relation
\begin{equation}\label{eq_Pythagoras}
||\mathbf{a}-\mathbf{c}||_2^2-||\mathbf{b}-\mathbf{c}||_2^2=||\mathbf{a}-\mathbf{b}||_2^2+2\langle\mathbf{a}-\mathbf{b},\mathbf{b}-\mathbf{c}\rangle,
\end{equation}
for the squared loss $h(\mathbf{x}_1,\cdots,\mathbf{x}_S)$. This property is much tighter than the property (\ref{eq_lipp}) of function with Lipschitz continuous gradient.

The result that using the squared loss leads to smaller Lipschitz constants by PIRE-PS is very interesting and useful. Intuitively, it results to minimize a tighter upper bounded surrogate function. Our experiments show that this will lead to a faster convergence of the PIRE-PS algorithm.

\subsection{PIRE with Alternative Updating}

In this section, we propose another splitting method to solve problem (\ref{eq_mspro}) based on the assumption that each  $\nabla_sh(\mathbf{x}_1,\cdots,\mathbf{x}_S)$ is Lipschitz continuous with constant $L_s(h)$. Different from PIRE-PS, which updates each $\mathbf{x}_s^{k+1}$ based on $\mathbf{x}_s^{k}$, $s=1,\cdots,S$, we instead update $\mathbf{x}_s^{k+1}$ based on all the latest $\mathbf{x}_s$. This is the known Gauss-Sidel iteration in numerical algebra. We name this method as PIRE with Alternative Updating (PIRE-AU).

Since $\nabla_sh(\mathbf{x}_1,\cdots,\mathbf{x}_S)$ is Lipschitz continuous, similar to (\ref{eq_lipp}), we have
\begin{equation}\label{eq_PIREnsmj}
\begin{split}
&h(\mathbf{x}_1^{k+1},\cdots,\mathbf{x}^{k+1}_{s-1},\mathbf{x}_s,\mathbf{x}_{s+1}^k,\cdots,\mathbf{x}_S^k)\\
\leq& h(\mathbf{x}_1^{k+1},\cdots,\mathbf{x}_{s-1}^{k+1},\mathbf{x}_{s}^{k},\cdots,\mathbf{x}_S^{k})+\\
&\langle\nabla_sh(\mathbf{x}_1^{k+1},\cdots,\mathbf{x}_{s-1}^{k+1},\mathbf{x}_s^k,\cdots,\mathbf{x}_S^k),\mathbf{x}_s-\mathbf{x}_s^k\rangle\\&+\frac{L_s(h)}{2}||\mathbf{x}_s-\mathbf{x}_s^k||_2^2.
\end{split}
\end{equation}
The hand right part of (\ref{eq_PIREnsmj}) is used as a surrogate function of $h(\mathbf{x}_1^{k+1},\cdots,\mathbf{x}^{k+1}_{s-1},\mathbf{x}_s,\mathbf{x}_{s+1}^k,\cdots,\mathbf{x}_S^k)$, which is tighter than (\ref{eq_lipp}) in PIRE. Then we update $\mathbf{x}_s^{k+1}$ by
\begin{equation}
\begin{split}\label{eq_updatexns}
\mathbf{x}_s^{k+1}=&\arg\min_{\mathbf{x_s}}\lambda\langle\mathbf{w}_s^k,\mathbf{g}_s(\mathbf{x}_s)\rangle+\frac{\mu_s}{2}||\mathbf{x}_s-\mathbf{x}_s^k||_2^2.\\
&+\langle\nabla_sh(\mathbf{x}_1^{k+1},\cdots,\mathbf{x}_{s-1}^{k+1},\mathbf{x}_s^k,\cdots,\mathbf{x}_S^k),\mathbf{x}_s-\mathbf{x}_s^k\rangle,
\end{split}
\end{equation}
where $\mu_s>L_s(h)/2$ and $\mathbf{w}_s^k$ is defined in (\ref{eq_updatewinpa}).

The updating rule in PIRE-AU by (\ref{eq_updatexns}) and (\ref{eq_updatewinpa}) also leads to converge. Any accumulation point of $\{\mathbf{x}^k\}$ is a stationary point. See the detailed proofs in the Supplementary Material.

Both PIRE-PS and PIRE-AU can solve the multi-variable problems. The advantage of PIRE-PS is that it is naturally parallelizable, while PIRE-AU may converge with less iterations due to smaller Lipschitz constants. If the squared loss function is used, PIRE-PS use the same small Lipschitz constants as PIRE-AU.




\begin{figure}
	\begin{subfigure}[b]{0.155\textwidth}
		\centering
        \includegraphics[width=\textwidth]{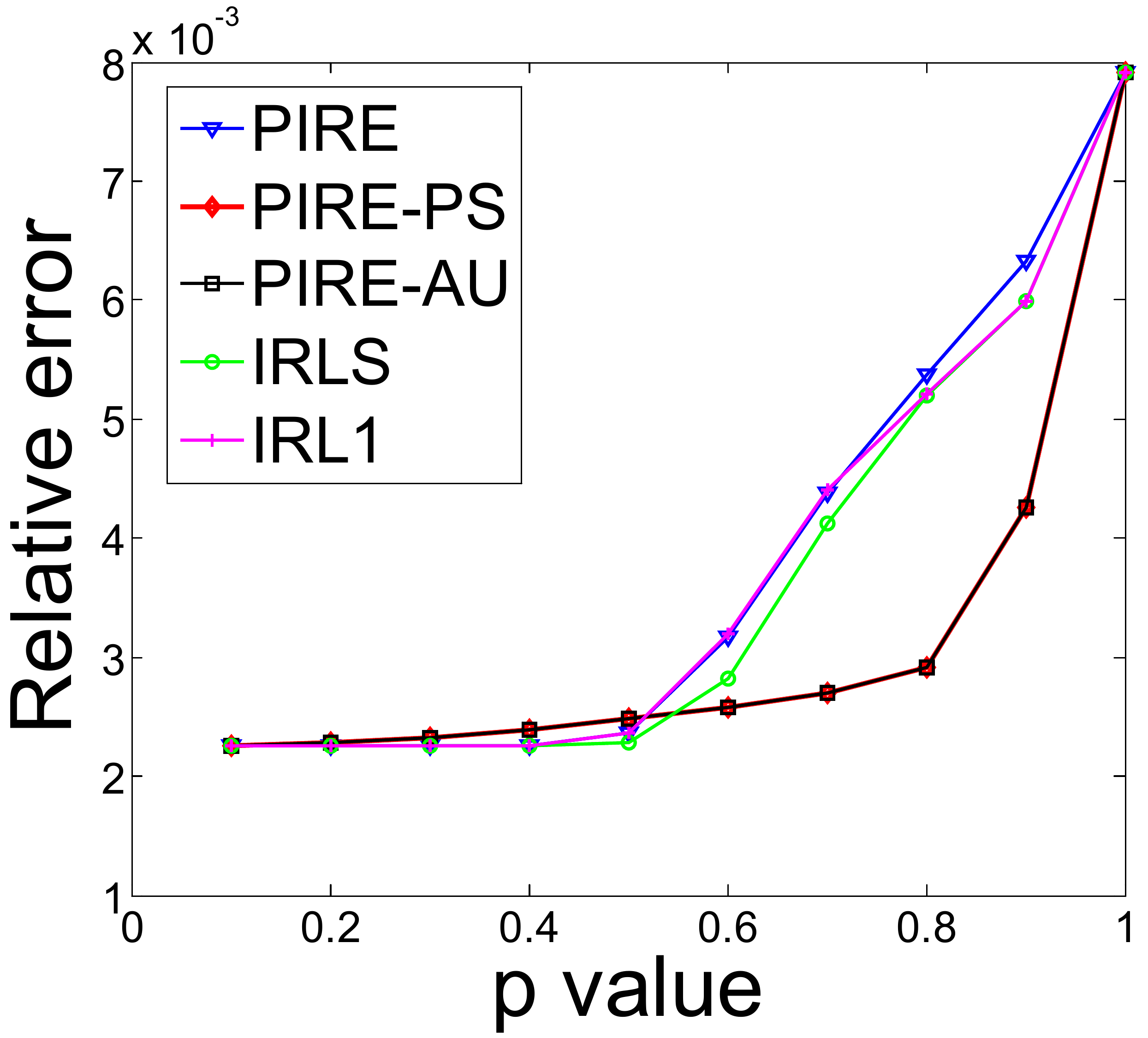}
        \caption{}
    \end{subfigure}
    \begin{subfigure}[b]{0.155\textwidth}
		\centering
        \includegraphics[width=\textwidth]{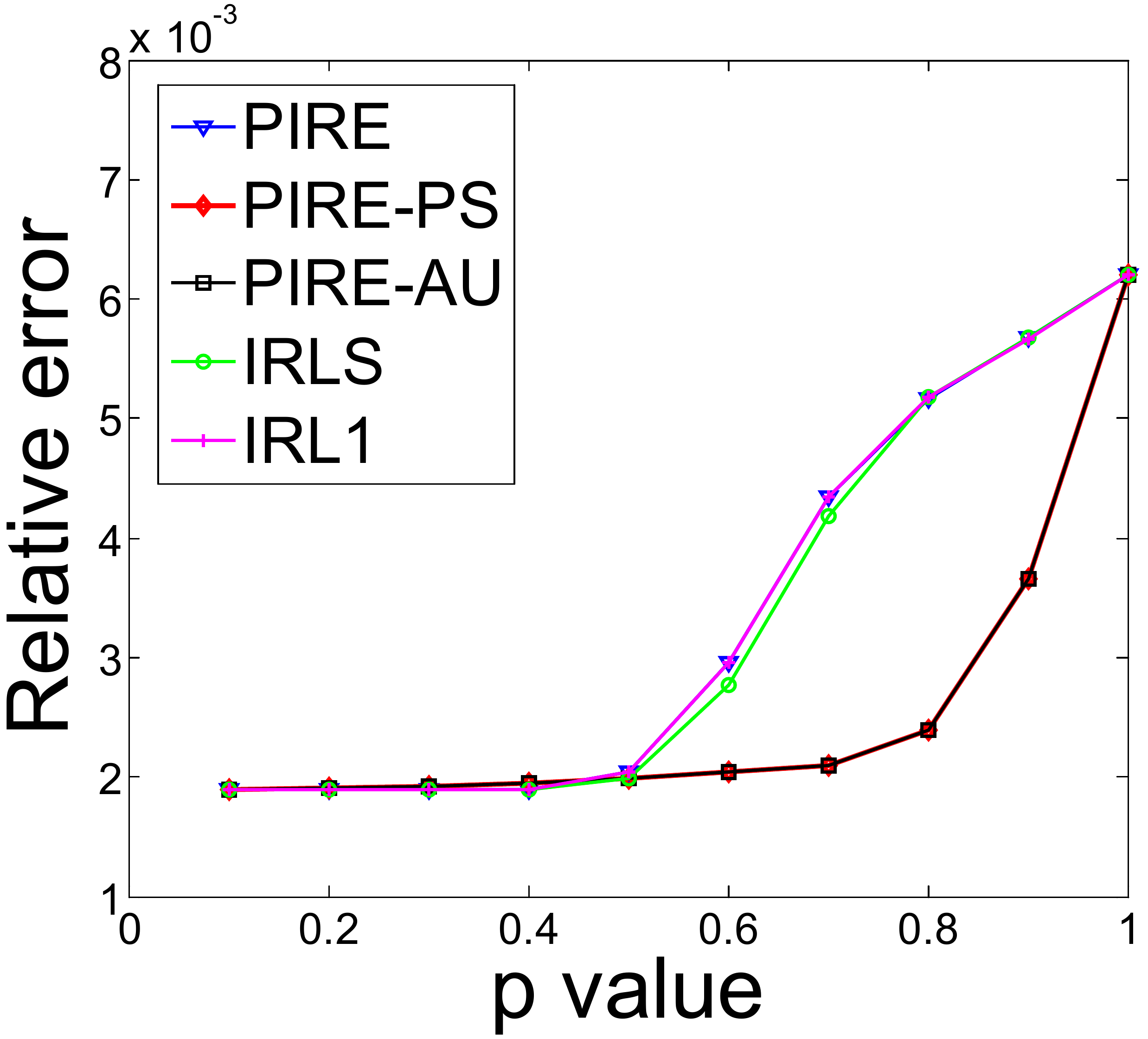}
        \caption{}
    \end{subfigure}
    \begin{subfigure}[b]{0.155\textwidth}
		\centering
        \includegraphics[width=\textwidth]{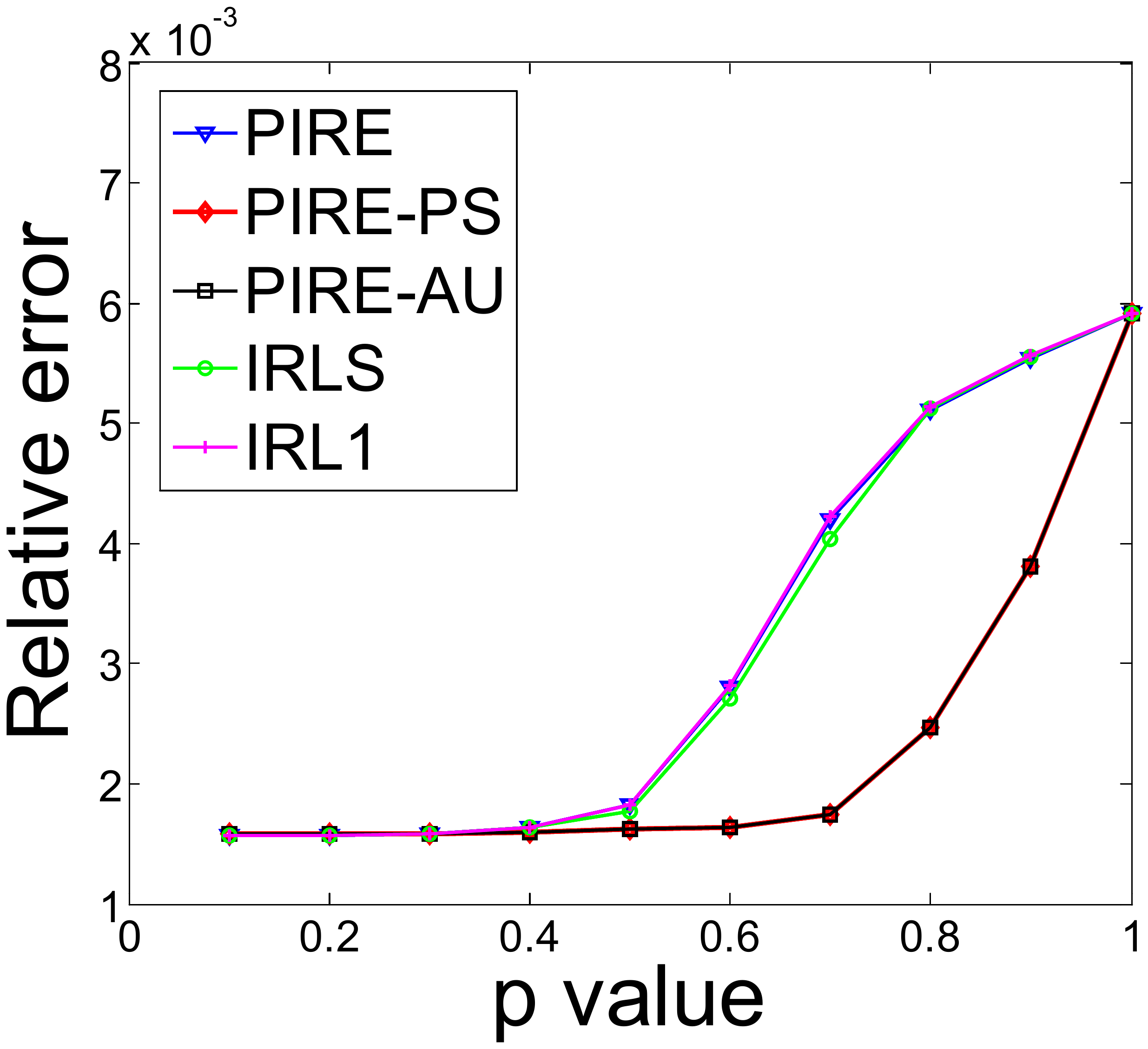}
        \caption{}
    \end{subfigure}
    \caption{\small Recovery performance comparison with a different number of measurement $\mathbf{A}\in\mathbb{R}^{m\times1000}$: (a) $m=200$; (b) $m=300$; and (c) $m=400$.}
\label{fig_recov_all_p}
\end{figure}

\section{Experiments}
\label{sec_experiments}
We present several numerical experiments to demonstrate the effectiveness of the proposed PIRE algorithm and its splitting versions. All the algorithms are implemented by Matlab, and are tested on a PC with 8 GB of RAM and Intel Core 2 Quad CPU Q9550.

\subsection{$\ell_p$-Minimization}
We compare our proposed PIRE, PIRE-PS and PIRE-AU algorithms with IRLS and IRL1 for solving the $\ell_p$-minimization problem (\ref{eq_lp}). For fair comparison, we try to use the same settings of all the completed algorithms. We use the solution to the $\ell_1$-minimization problem as the initialization. We find that this will accelerate the convergence of the iteratively reweighted algorithms, and also enhance the recovery performance. The choice of $\epsilon$ in (\ref{eq_smo11}) and (\ref{eq_IRLSp}) plays an important role for sparse signal recovery, but theoretical support has not been carried out so far. Several different decreasing rules have been tested before \cite{candes2008enhancing,IRLSrank,lai2013improved}, but none of them dominates others. Since the sparsity of sparse signal is usually unknown, we empirically set $\epsilon^{k+1}=\epsilon^k/\rho$, with $\epsilon^0=0.01$, and $\rho=1.1$ \cite{IRLSrank}. The algorithms are stopped when $||\mathbf{x}^k-\mathbf{x}^{k+1}||_2/||\mathbf{x}^k||_2\leq10^{-6}$.

IRL1 requires solving (\ref{eq_IRL1w}) as inner loop. FISTA is employed to solve (\ref{eq_IRL1w}) with warm start, i.e. using $\mathbf{x}^k$ as initialization to obtain $\mathbf{x}^{k+1}$. This trick greatly reduces the inner loop iteration, which is the main cost for IRL1. For PIRE-PS and PIRE-AU algorithms, we solve problem (\ref{eq_mspro}) by setting $S=20$.

\begin{table}[!t]
\tiny
\centering
\caption{\small Comparison of iteration number, running time (in seconds), objective function value and relative recovery error of different iterative reweighted methods.}
\label{Tab_timecom}
\centering
\begin{tabular}{c| l| l| r@{.}l| c| c }
\hline
\multirow{2}*{ Size ($m,n,t$) }
 & \multirow{2}*{ Methods }  & \multirow{2}*{ Iter. }  & \multicolumn{2}{|c|}{Time} & Obj.  & Recovery error   \\
 & & & (second)&&  ($\times 10^{-2}$) & ($\times10^{-3}$) \\\hline
\multirow{5}*{ (100,500,50) }
 &       PIRE &    116  &  0 & 70 &  \textbf{5.238} &  2.529   \\
 &    PIRE-PS &     58  &  \textbf{0} & \textbf{48} &  5.239 &  2.632   \\
 &    PIRE-AU &     \textbf{56}  &  0 & 63 &  5.239 &  2.632   \\
 &       IRLS &    168  & 81 & 82 &  5.506 &  \textbf{2.393}   \\
 &       IRL1 &     \textbf{56}  &  3 & 43 &  5.239 &  2.546   \\\cline{1-7}
\multirow{5}*{ (200,800,100) }
 &       PIRE &    119  &  1 & 48 & 16.923 &  2.246   \\
 &    PIRE-PS &     37  &  \textbf{0} & \textbf{82} & \textbf{16.919} &  2.192   \\
 &    PIRE-AU &     \textbf{36}  &  0 & 88 & \textbf{16.919} &  2.192   \\
 &       IRLS &    169  & 474 & 19 & 17.784 &  \textbf{2.142}   \\
 &       IRL1 &     81  & 13 & 53 & 16.924 &  2.248 \\\cline{1-7}
\multirow{5}*{ (300,1000,200) }\
 &       PIRE &    151  &  4 & 63 & 42.840 &  2.118   \\
 &    PIRE-PS &     29  &  1 & 38 & \textbf{42.815} &  1.978   \\
 &    PIRE-AU &     \textbf{28}  &  \textbf{1} & \textbf{34} & \textbf{42.815} &  \textbf{1.977}   \\
 &       IRLS &    171  & 1298 & 70 & 44.937 &  2.015   \\
 &       IRL1 &     79  & 35 & 59 & 42.844 &  2.124  \\\cline{1-7}
\multirow{5}*{ (500,1500,200) }
 &       PIRE &    159  &  8 & 88 & 64.769 &  2.010   \\
 &    PIRE-PS &     26  &  2 & 27 & \textbf{64.718} &  \textbf{1.814}   \\
 &    PIRE-AU &     \textbf{25}  &  \textbf{2} & \textbf{20} & \textbf{64.718} &  \textbf{1.814}   \\
 &       IRLS &    171  & 3451 & 79 & 67.996 &  1.923   \\
 &       IRL1 &     89  & 80 & 89 & 64.772 &  2.013    \\\cline{1-7}
\multirow{5}*{ ( 800,2000,200) }
 &       PIRE &    140  & 14 & 99 & 87.616 &  1.894   \\
 &    PIRE-PS &     33  &  5 & 15 & \textbf{87.533} &  \textbf{1.648}   \\
 &    PIRE-AU &     \textbf{32}  &  \textbf{4} & \textbf{97} & \textbf{87.533} &  \textbf{1.648}   \\
 &       IRLS &    177  & 7211 &  2 & 91.251 &  1.851   \\
 &       IRL1 &    112  & 173 & 26 & 87.617 &  1.895  \\\cline{1-7}
\end{tabular}
\end{table}

\subsubsection{Sparse Signal Recovery}
The first experiment is to examine the recovery performance of sparse signals by using the proposed methods. The setup for each trial is as follows. The dictionary $\mathbf{A}\in\mathbb{R}^{m\times n}$ is a Gaussian random matrix generated by Matlab command \mcode{randn}, with the sizes $m=200,300,400$, and $n=1000$. The sparse signal $\mathbf{x}$ is randomly generated with sparsity $||\mathbf{x}||_0=20$. The response $\mathbf{b}=\mathbf{A}\mathbf{x}+0.01\mathbf{e}$, where $\mathbf{e}$ is Gaussian random vector. Given $\mathbf{A}$ and $\mathbf{b}$, we can recover $\hat{\mathbf{x}}$ by solving the $\ell_p$-minimization problem by different methods. The parameter is set to $\lambda=10^{-4}$. We use the relative recovery error $||\hat{\mathbf{x}}-\mathbf{x}||_2/||\mathbf{x}||_2$ to measure the recovery performance. Based on the above settings and generated data, we find that the recovery performances are stable. We run 20 trials and report the mean relative error for comparison.

Figure \ref{fig_recov_all_p} plots the relative recovery errors v.s. different $p$ values ($p=0.1,\cdots,0.9,1$) on three data sets with different numbers of measurements. The result for $p=1$ is obtained by FISTA for $\ell_1$-minimization. We can see that all the iteratively reweighted algorithms achieve better recovery performance with $p<1$ than $\ell_1$-minimization. Also a smaller value of $p$ leads to better recovery performance, though the $\ell_p$-minimization problem is nonconvex and a globally optimal solution is not available. In most cases, PIRE is comparative with IRLS and IRL1. A surprising result is that   PIRE-PS and PIRE-AU outperform the other methods when $0.5<p<1$. They use a smaller Lipschitz constant than PIRE, and thus may converge faster. But none of these iteratively reweighted methods is guaranteed to be optimal.


\begin{figure}
	\begin{subfigure}[b]{0.155\textwidth}
		\centering
        \includegraphics[width=\textwidth]{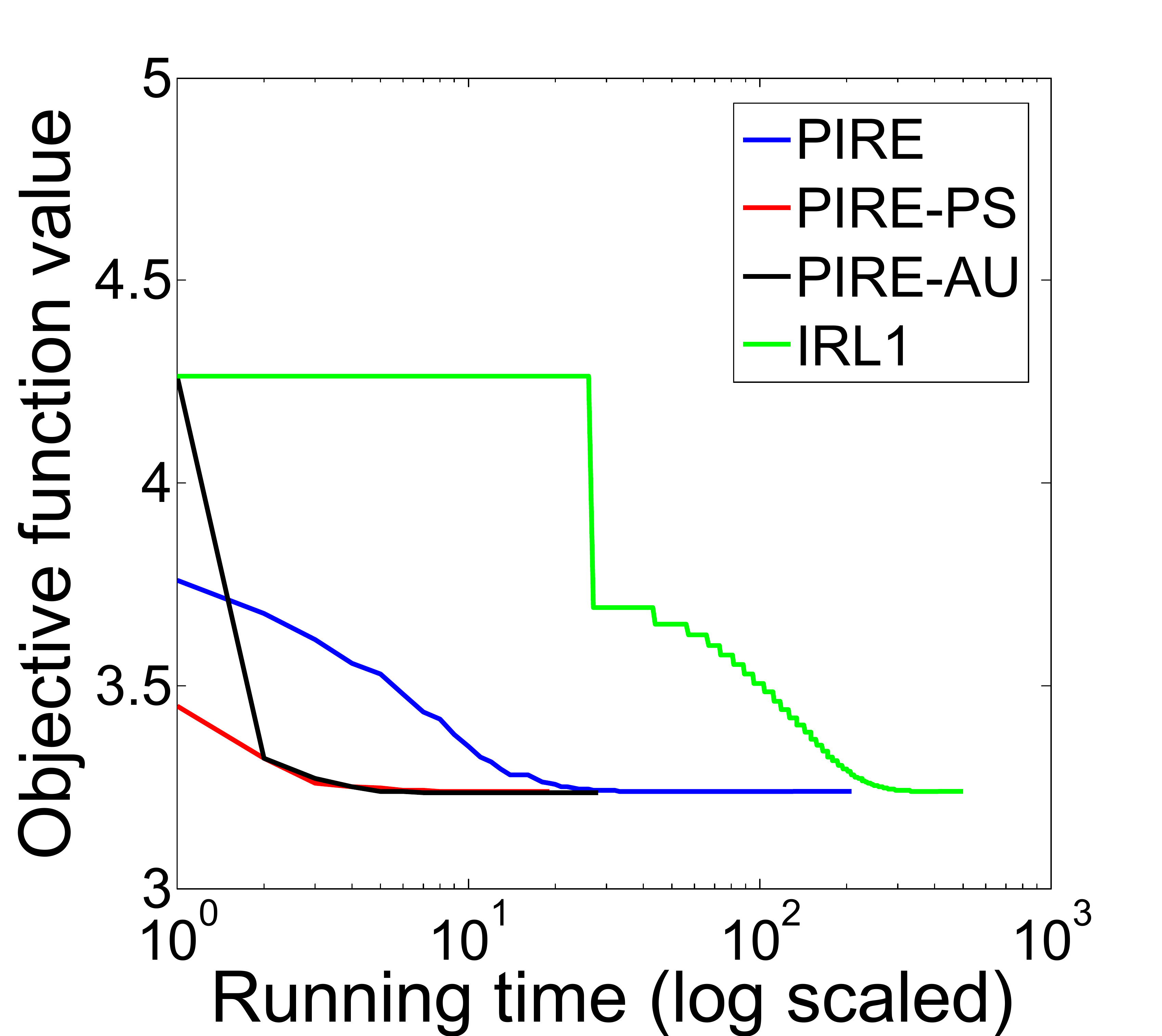}
        \caption{}
    \end{subfigure}
    \begin{subfigure}[b]{0.155\textwidth}
		\centering
        \includegraphics[width=\textwidth]{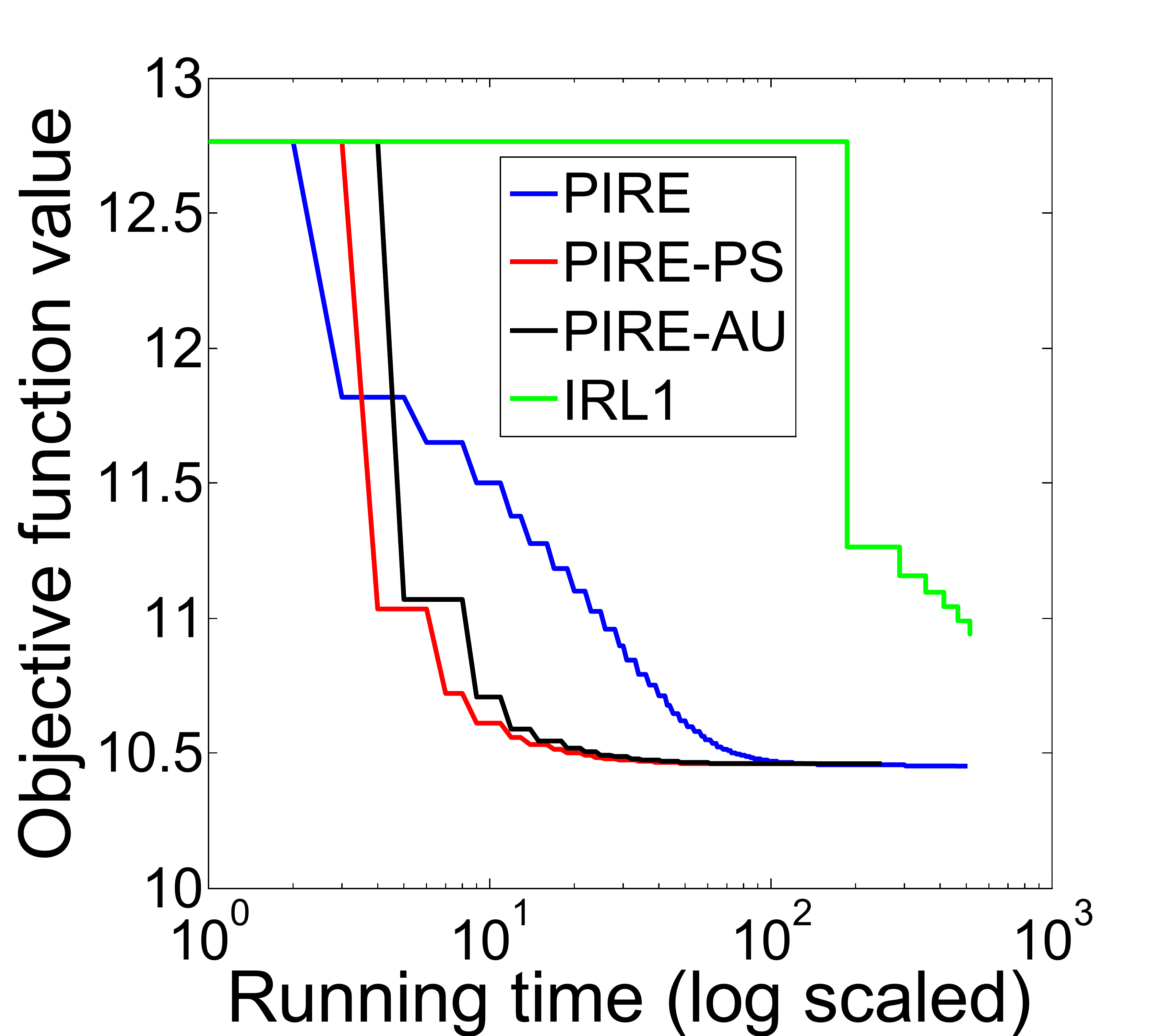}
        \caption{}
    \end{subfigure}
    \begin{subfigure}[b]{0.155\textwidth}
		\centering
        \includegraphics[width=\textwidth]{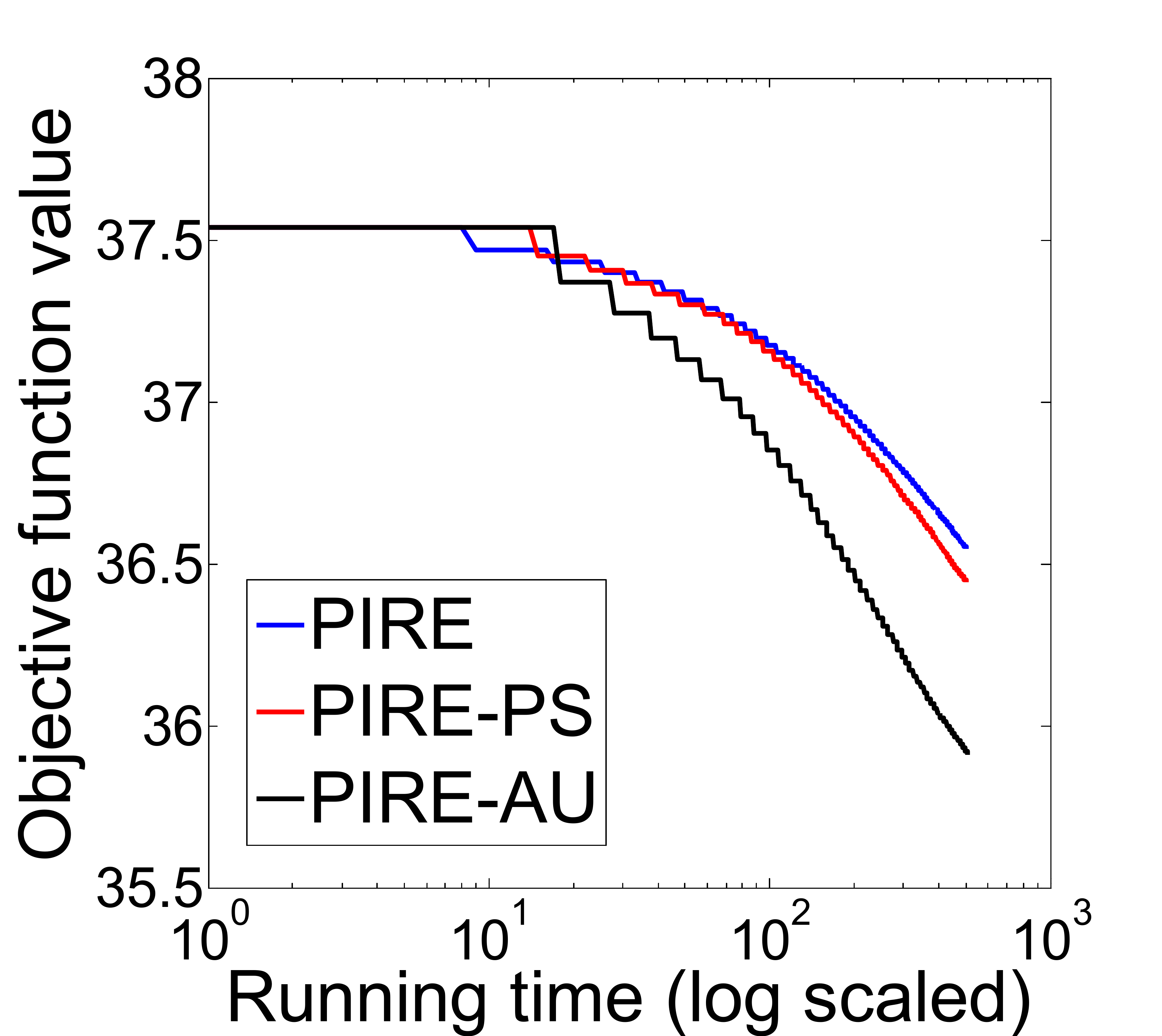}
        \caption{}
    \end{subfigure}
    \caption{\small Running time v.s. objective function value on three synthesis data sets with size $(m,n,t)$: (a) (1000,3000,500); (b) (1000,5000,1000); (c) (1000,10000,1000).}
\label{fig_timobj}
\vspace{-1em}
\end{figure}

\subsubsection{Running Time Comparison}
The second experiment is to show the advantage in running time of the proposed methods. We implement all the completed methods in
matrix form for solving the following $\ell_p$-minimization problem
\begin{equation}
\min_{\mathbf{X}\in\mathbb{R}^{n\times t}} \lambda||\mathbf{X}||_p^p+\frac{1}{2}||\mathbf{A}\mathbf{X}-\mathbf{B}||_F^2,
\end{equation}
where $\mathbf{A}\in\mathbb{R}^{m\times n}$ and $\mathbf{B}\in\mathbb{R}^{m\times t}$, $||\mathbf{X}||_p^p=\sum_{ij}|X_{ij}|^p$, and $p$ is set to 0.5 in this test. $\mathbf{A}$, $\mathbf{B}$ and $\mathbf{X}$ are generated by the same procedure as the above section, and the same settings of algorithm parameters are followed. Each column of $\mathbf{X}$ is with sparsity $n\times2\%$. We test on several different sizes of data sets, parameterized as $(m,n,t)$. The iteration number, running time, objective function value and the relative recovery error are tabulated in Table \ref{Tab_timecom}. It can be seen that the proposed methods are much more efficient than IRLS and IRL1. PIRE-PS and PIRE-AU converge with less iteration and less running time. In our test, IRL1 is more efficient than IRLS. The reasons lie in: (1) initialization as a sparse solution to $\ell_1$-minimization is a good choice for IRL1, but not for IRLS; (2) For each iteration in IRLS, solving $t$ equations (\ref{eq_sol2}) in a loop by Matlab is not efficient; (3) IRL1 converges with less inner loop iterations due to warm start.

We also plot the running time v.s. objective function value on three larger data sets in Figure \ref{fig_timobj}. The algorithms are stopped within 500 seconds in this test. IRLS costs much more time, and thus it is not plotted. IRL1 is not plotted for the case $n=10,000$. It can be seen that PIRE-PS and PIRE-AU decreases the objective function value faster than PIRE.

\subsection{Multi-Task Feature Learning}

In this experiment, we use our methods to solve the multi-task learning problem. Assume we are given $m$ learning tasks associated with $\{(\mathbf{X}_1,\mathbf{y}_1),\cdots,(\mathbf{X}_m,\mathbf{y}_m)\}$, where $\mathbf{X}_i\in\mathbb{R}^{n_i\times d}$ is the data matrix of the $i$-th task with each row a sample, $\mathbf{y}_i\in\mathbb{R}^{n_i}$ is the label of the $i$-th task, $n_i$ is the number of samples for the $i$-th task, and $d$ is the data dimension. Our goal is to find a matrix $\mathbf{Z}=[\mathbf{z}_1,\cdots,\mathbf{z}_m]\in\mathbb{R}^{d\times m}$ such that $\mathbf{y}_i\approx\mathbf{X}_i\mathbf{z}_i$.  The capped-$\ell_1$ norm is used to regularize $\mathbf{Z}$ \cite{gong2012multi}
\begin{equation}\label{multitask}
\min_{\mathbf{Z}}\lambda\sum_{j=1}^d\min(||\mathbf{z}^j||_1,\theta)+h(\mathbf{Z}),
\end{equation}
where $h(\mathbf{Z})=\sum_{i=1}^m||\mathbf{X}_i\mathbf{z}_i-\mathbf{y}_i||_2^2/mn_i$ is the loss function, $\theta>0$ is the thresholding parameter, and $\mathbf{z}^j$ is the $j$-th row of $\mathbf{Z}$. The above problem can be solved by our proposed PIRE, PIRE-PS and PIRE-AU algorithms, by letting $f(\mathbf{y})=\sum_{j=1}^d\min(y_i,\theta)$, and $g(\mathbf{Z})=[||\mathbf{z}^1||_1;\cdots;||\mathbf{z}^m||_1]$.

\begin{figure}
	\begin{subfigure}[b]{0.23\textwidth}
		\centering
        \includegraphics[width=\textwidth]{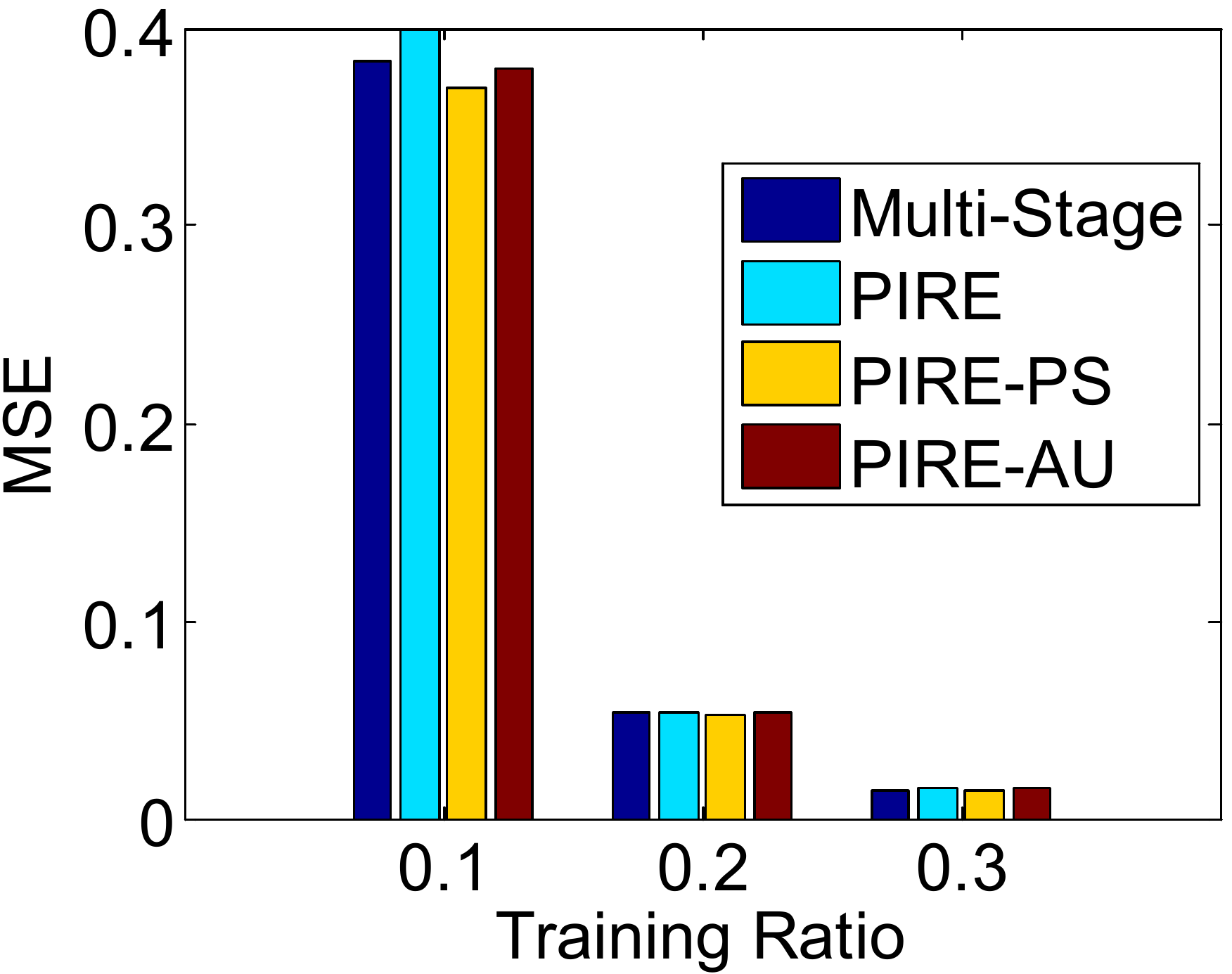}
        \caption{}
    \end{subfigure}
    \begin{subfigure}[b]{0.23\textwidth}
		\centering
        \includegraphics[width=\textwidth]{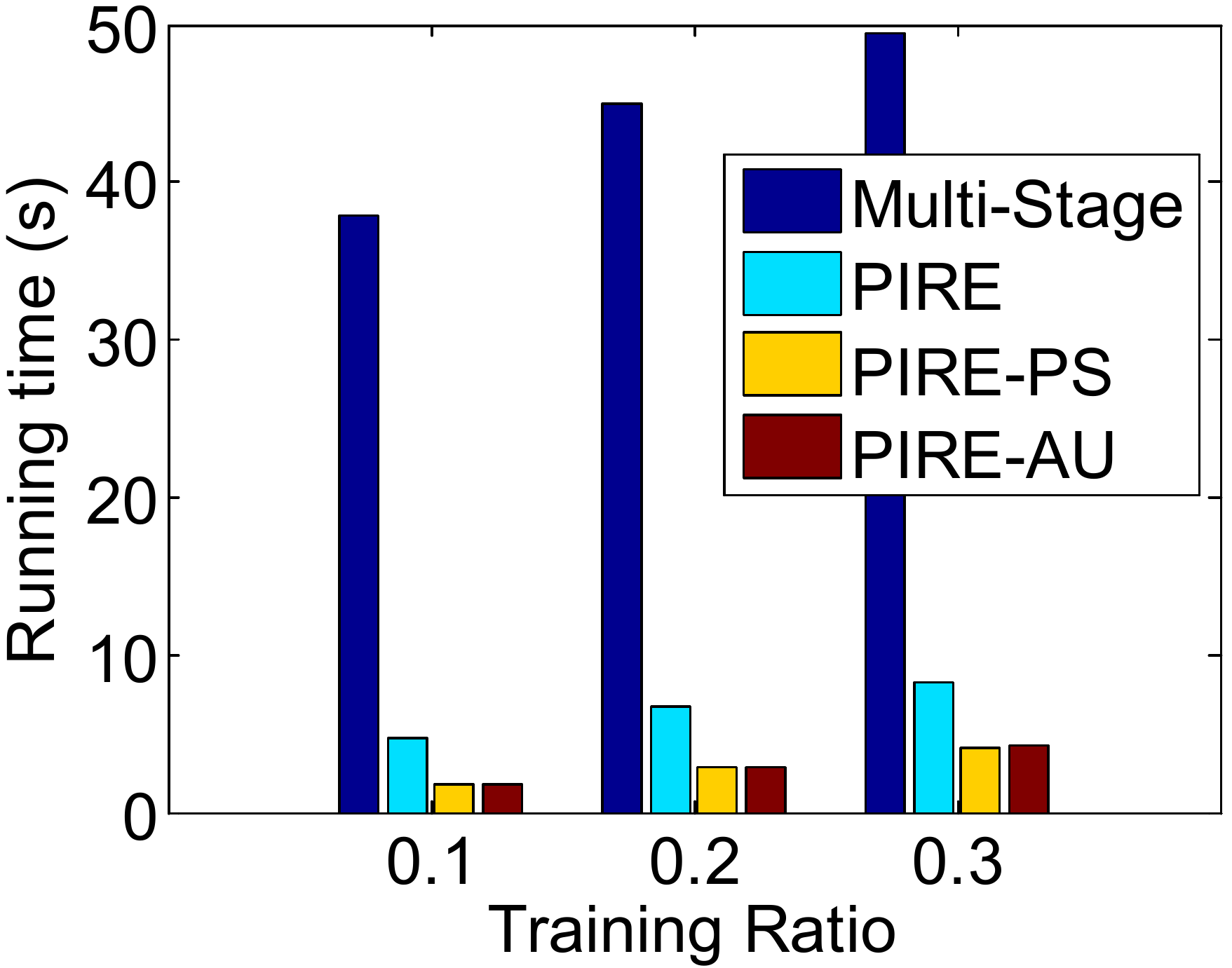}
        \caption{}
    \end{subfigure}
    \caption{\small Comparison of (a) mean squared error (MSE) and running time on the Isolet data set for multi-task feature learning.}
\vspace{-1em}
\label{fig_res_Isolet}
\end{figure}

The Isolet \cite{UCI} data set is used in our test. 150 subjects spoke the name of each letter of the alphabet twice. Hence, we have 52 training examples from each speaker. The speakers are grouped into 5 subsets of 30 speakers each. Thus, we have 5 tasks with each task corresponding to a subset. There are 1560, 1560, 1560, 1558, and 1559 samples of 5 tasks, respectively. The data dimension is 617, and the response is the English letter label (1-26). We randomly select the training samples from each task with different training ratios (0.1, 0.2 and 0.3) and use the rest of samples to form the test set. We compare our PIRE, PIRE-PS and PIRE-AU (we set $S=m=5$ in PIRE-PS and PIRE-AU) with the Multi-Stage algorithm \cite{zhang2008multi}. We report the Mean Squared Error (MSE) on the test set and the running time for solving (\ref{multitask}) on the training set. The results are averaged over 10 random splittings. As shown in Figure \ref{fig_res_Isolet}, it can be seen that all these methods achieve comparative performance, but our PIRE, PIRE-PS and PIRE-AU are much more efficient than the Multi-Stage algorithm.

\section{Conclusions}
This paper proposes the PIRE algorithm for solving the general problem (\ref{eq_genp}). PIRE solves a series of problem (\ref{eq_proxw}), whose computational cost is usually very cheap. We further propose two splitting versions of PIRE to handle the multi-variable problems. In theory, we prove that  PIRE (also its splitting versions) converges and any limit point is a stationary point. We test our methods to solve the $\ell_p$-minimization problem and multi-task feature learning problem. Experimental results on both synthesis and real data sets show that our methods are with comparative learning performance, but much more efficient, by comparing with IRLS and IRL1 or multi-stage algorithms.
It would be interesting to apply PIRE for structured sparsity optimization, and also the nonconvex low rank regularized minimization problems \cite{IRNN}.

\section*{Acknowledgements}

This research is supported by the Singapore National Research Foundation under its International Research Centre @Singapore Funding Initiative and administered by the IDM Programme Office. Z. Lin is supported by NSF of China (Grant nos. 61272341, 61231002, and 61121002) and MSRA.

\bibliography{PIRE}
\bibliographystyle{aaai}

\newpage
\onecolumn
~\\\\
\begin{center}
\LARGE\textbf{Supplementary Material of \\Proximal Iteratively Reweighted Algorithm with Multiple Splitting for Nonconvex Sparsity Optimization}
\end{center}

\vspace{1em}
\renewcommand{\thefootnote}{$\star$} 
\begin{center}
{\Large\textbf{Canyi Lu}$^1$, \textbf{Yunchao Wei}$^2$, \textbf{Zhouchen Lin}$^{3,}$\footnote{Corresponding author.}, \textbf{Shuicheng Yan}$^1$}\\
$^1$ Department of Electrical and Computer Engineering, National University of Singapore\\
$^2$ Institute of Information Science, Beijing Jiaotong University\\
$^3$ Key Laboratory of Machine Perception (MOE), School of EECS, Peking University\\
{\tt\small canyilu@gmail.com, wychao1987@gmail.com, zlin@pku.edu.cn, eleyans@nus.edu.sg}
\end{center}

\vspace{0.01em}
\vspace{1em}

In this supplementary material, we give the detailed convergence results and proofs of the PIRE-PS and PIRE-AU algorithms.

\section{Proof of Theorem \ref{thm_con}}
\textbf{Proof}. The sequence $\{\mathbf{x}^k\}$ generated in Algorithm \ref{alg_PIRE} is bounded by Theorem \ref{thm_pro}. Hence, there exists an accumulation point $\mathbf{x}^*$, and a subsequence $\{\mathbf{x}^{k_{j}}\}$ such that $\lim\limits_{j\rightarrow\infty}\mathbf{x}^{k_{j}}=\mathbf{x}^*$. From the fact that $\lim\limits_{k\rightarrow\infty}(\mathbf{x}^{k}-\mathbf{x}^{k+1})=\bm{0}$ in Theorem \ref{thm_pro}, we have $\lim\limits_{j\rightarrow\infty}\mathbf{x}^{k_{j}+1}=\mathbf{x}^*$. Since $\mathbf{x}^{k_j+1}$ solves problem (\ref{eq_updatex}), there exists $\mathbf{v}^{k_j+1}\in\partial\langle\mathbf{w}^{k_j},\mathbf{g}(\mathbf{x}^{k_j+1})\rangle$ such that
\begin{equation}\label{eq_profthm2}
\lambda\mathbf{v}^{k_j+1}+\nabla h(\mathbf{x}^{k_j})+\mu(\mathbf{x}^{k_j+1}-\mathbf{x}^{k_j})=\mathbf{0}.
\end{equation}
By the upper semi-continuous property of the subdifferential \cite{clarke1983nonsmooth}, there exists $\mathbf{v}^*\in\partial\langle\mathbf{w}^{*},\mathbf{g}(\mathbf{x}^{*})\rangle$, with $\mathbf{w}^*\in-\partial\left( -f(\mathbf{g}(\mathbf{x}^*))\right)$, such that
\begin{equation}
\mathbf{0}=\lambda\mathbf{v}^*+\nabla h(\mathbf{x}^*)\in\frac{\partial F(\mathbf{x}^*)}{\partial \mathbf{x}^*}.
\end{equation}
Hence $\mathbf{x}^*$ is a stationary point of (\ref{eq_genp}).

Furthermore, summing (\ref{eq_proof66}) for $k=1,\cdots,n$, we have
\begin{equation}
\begin{split}
F(\mathbf{x}^1)-F(\mathbf{x}^{n+1})&\geq\left(\mu-\frac{L(h)}{2}\right)\sum_{k=1}^n||\mathbf{x}^{k+1}-\mathbf{x}^k||_2^2\\
&\geq n\left(\mu-\frac{L(h)}{2}\right)\min_{1\leq k\leq n}||\mathbf{x}^{k+1}-\mathbf{x}^k||_2^2.
\end{split}
\end{equation}
Thus
\begin{equation}
\begin{split}
\min_{1\leq k\leq n}||\mathbf{x}^{k+1}-\mathbf{x}^k||_2^2\leq&\frac{F(\mathbf{x}^1)-F(\mathbf{x}^{n+1})}{n\left(\mu-\frac{L(h)}{2}\right)}\\
\leq&\frac{F(\mathbf{x}^1)-F(\mathbf{x}^{*})}{n\left(\mu-\frac{L(h)}{2}\right)}.
\end{split}
\end{equation}
$\hfill\blacksquare$

\section{Convergence Analysis of PIRE-PS}
For the general PIRE-PS algorithm for solving problem (\ref{eq_mspro}) by (\ref{updatexinps}) and (\ref{eq_updatewinpa}), actually its updating rule is the same as PIRE. Thus the same convergence results in Theorem \ref{thm_pro} and \ref{thm_con} hold. In this section, we provide the convergence analysis of PIRE-PS algorithm by (\ref{eq_psupx}) and (\ref{eq_updatewinpa}) for solving problem (\ref{eq_mspro}) with the squared loss function $h(\mathbf{x}_1,\cdots,\mathbf{x}_S)=\frac{1}{2}\left\|\sum_{s=1}^S\mathbf{A}_s\mathbf{x}_s-\mathbf{b}\right\|_2^2$.

\begin{theo}\label{thm_ps_pro_ps}
Assume that the squared loss function $h(\mathbf{x}_1,\cdots,\mathbf{x}_S)=\frac{1}{2}\left\|\sum_{s=1}^S\mathbf{A}_s\mathbf{x}_s-\mathbf{b}\right\|_2^2$ is used in problem (\ref{eq_mspro}).
Let $D=F(\mathbf{x}^1)$, and $\mu_s>\frac{||\mathbf{A}_s||_2^2}{2}$, $s=1,\cdots,S$, and $\delta=\min_s\{\mu_s-||\mathbf{A}_s||_2^2/2\}$. The sequence $\{\mathbf{x}^k\}$ generated by PIRE-PS in (\ref{eq_psupx}) and (\ref{eq_updatewinpa}) satisfies the following properties:
\begin{enumerate}[(1)]
\item $F(\mathbf{x}^k)$ is monotonically decreasing, i.e. $F(\mathbf{x}^{k+1})\leq F(\mathbf{x}^k)$. Indeed,
\begin{equation*}
F(\mathbf{x}^k)-F(\mathbf{x}^{k+1})\geq\sum_{s=1}^S\left(\mu_s-\frac{||\mathbf{A}_s||_2^2}{2}\right)||\mathbf{x}_s^{k+1}-\mathbf{x}_s^k||^2\geq0;
\end{equation*}
\item The sequence $\{\mathbf{x}^k\}$ is bounded;
\item $\lim\limits_{k\rightarrow\infty}(\mathbf{x}^k-\mathbf{x}^{k+1})=\bm{0}$.
\item Any accumulation point of $\{\mathbf{x}^k\}$ is a stationary point $\mathbf{x}^*$ of problem (\ref{eq_mspro}).
\item For any $n\geq1$, $\min_{1\leq k\leq n}||\mathbf{x}^{k+1}-\mathbf{x}^{k}||_2^2\leq\frac{F(\mathbf{x}^1)-F(\mathbf{x}^*)}{n\delta}$.
\end{enumerate}
\end{theo}
\textbf{Proof}. Since $\mathbf{x}^{k+1}_s$ is the globally optimal solution to problem (\ref{eq_psupx}), the zero vector is contained in the subgradient with respect to $\mathbf{x}_s$. That is, there exists $\mathbf{v}_s^{k+1}\in\partial\langle\mathbf{w}_s^k,\mathbf{g}(\mathbf{x}^{k+1}_s)\rangle$ such that
\begin{equation}\label{eq_proof1_ps111}
\lambda\mathbf{v}^{k+1}_s+\mathbf{A}_s^T(\mathbf{A}\mathbf{x}^k-\mathbf{b})+\mu_s(\mathbf{x}^{k+1}_s-\mathbf{x}^{k}_s)=\mathbf{0}.
\end{equation}
A dot-product with $\mathbf{x}^{k+1}_s-\mathbf{x}^k_s$ on both sides of (\ref{eq_proof1_ps111}) gives
\begin{equation}\label{eq_proof2_ps222}
\begin{split}
\lambda\left\langle \mathbf{v}_s^{k+1},\mathbf{x}_s^{k+1}-\mathbf{x}_s^k\right\rangle+\left\langle \mathbf{A}\mathbf{x}^k-\mathbf{b},\mathbf{A}_s(\mathbf{x}_s^{k+1}-\mathbf{x}_s^k)\right\rangle+\mu_s||\mathbf{x}_s^{k+1}-\mathbf{x}_s^k||_2^2=0.
\end{split}
\end{equation}
Recalling the definition of the subgradient of the convex function, we have
\begin{equation}\label{eq_proof3_ps}
\begin{split}
\langle\mathbf{w}^k_s,\mathbf{g}(\mathbf{x}_s^k)-\mathbf{g}(\mathbf{x}_s^{k+1})\rangle\geq\left\langle\mathbf{v}^{k+1}_s,\mathbf{x}_s^k-\mathbf{x}_s^{k+1}\right\rangle.
\end{split}
\end{equation}
Combining (\ref{eq_proof2_ps222}) and (\ref{eq_proof3_ps}) gives
\begin{equation}\label{eq_psproof1_ps}
\begin{split}
&\lambda\sum_{s=1}^S\langle\mathbf{w}_s^k,\mathbf{g}(\mathbf{x}_s^k)-\mathbf{g}(\mathbf{x}_s^{k+1})\rangle\\
\geq&-\sum_{s=1}^S\left\langle\mathbf{A}\mathbf{x}^k-\mathbf{b},\mathbf{A}_s(\mathbf{x}_s^k-\mathbf{x}_s^{k+1})\right\rangle+\sum_{s=1}^S\mu_s||\mathbf{x}_s^{k+1}-\mathbf{x}_s^k||_2^2\\
=&-\langle\mathbf{A}\mathbf{x}^k-\mathbf{b},\mathbf{A}(\mathbf{x}^k-\mathbf{x}^{k+1})\rangle+\sum_{s=1}^S\mu_s||\mathbf{x}_s^{k+1}-\mathbf{x}_s^k||_2^2.
\end{split}
\end{equation}
By using the Pythagoras relation
\begin{equation}\label{eq_Pythagoras_ps}
||\mathbf{a}-\mathbf{c}||^2_2-||\mathbf{b}-\mathbf{c}||^2_2=||\mathbf{a}-\mathbf{b}||^2_2+2\langle\mathbf{a}-\mathbf{b},\mathbf{b}-\mathbf{c}\rangle,
\end{equation}
we get
\begin{equation}\label{eq_psproof2_ps}
\begin{split}
&\frac{1}{2}\left(\left\|\mathbf{A}\mathbf{x}^k-\mathbf{b}\right\|_2^2-\left\|\mathbf{A}\mathbf{x}^{k+1}-\mathbf{b}\right\|_2^2\right)
=\frac{1}{2}\left\|\mathbf{A}\left(\mathbf{x}^k-\mathbf{x}^{k+1}\right)\right\|^2_2+\left\langle\mathbf{A}(\mathbf{x}^k-\mathbf{x}^{k+1}),\mathbf{A}\mathbf{x}^{k+1}-\mathbf{b}\right\rangle.\\
\end{split}
\end{equation}
By the assumption (\textbf{C1}) that $f_s$ is concave, we have
\begin{equation}\label{eq_like22_ps}
f_s(\mathbf{g}_s(\mathbf{x}_s^{k}))-f_s(\mathbf{g}_s(\mathbf{x}_s^{k+1}))\geq\langle\mathbf{w}_s^k,\mathbf{g}_s(\mathbf{x}_s^{k})-\mathbf{g}_s(\mathbf{x}_s^{k+1})\rangle.
\end{equation}
Now, combining (\ref{eq_psproof1_ps})(\ref{eq_psproof2_ps}) and (\ref{eq_like22_ps}) leads to
\begin{equation}\label{eq_proofPIREPS111}
\begin{split}
&F(\mathbf{x}^k)-F(\mathbf{x}^{k+1})\\
=&\lambda\sum_{s=1}^S\left(f_s(\mathbf{g}_s(\mathbf{x}_s^{k}))-f_s(\mathbf{g}_s(\mathbf{x}_s^{k+1}))\right)+\frac{1}{2}\left(\left\|\mathbf{A}\mathbf{x}^k-\mathbf{b}\right\|_2^2-\left\|\mathbf{A}\mathbf{x}^{k+1}-\mathbf{b}\right\|_2^2\right)\\
\geq&-\langle\mathbf{A}\mathbf{x}^k-\mathbf{b},\mathbf{A}(\mathbf{x}^k-\mathbf{x}^{k+1})\rangle+\sum_{s=1}^S\mu_s||\mathbf{x}_s^{k+1}-\mathbf{x}_s^k||_2^2+\frac{1}{2}\left\|\mathbf{A}\left(\mathbf{x}^k-\mathbf{x}^{k+1}\right)\right\|_2^2+\left\langle\mathbf{A}(\mathbf{x}^k-\mathbf{x}^{k+1}),\mathbf{A}\mathbf{x}^{k+1}-\mathbf{b}\right\rangle\\
=&-\frac{1}{2}\left\|\mathbf{A}\left(\mathbf{x}^k-\mathbf{x}^{k+1}\right)\right\|_2^2+\sum_{s=1}^S\mu_s||\mathbf{x}_s^{k+1}-\mathbf{x}_s^k||_2^2\\
\geq&\sum_{s=1}^S\left(\mu_s-\frac{||\mathbf{A}_s||_2^2}{2}\right)||\mathbf{x}_s^{k+1}-\mathbf{x}_s^k||_2^2\\
\geq&\delta||\mathbf{x}^{k+1}-\mathbf{x}^k||_2^2\geq0,
\end{split}
\end{equation}
Thus $F(\mathbf{x}^k)$ is monotonically decreasing. Summing all the
above inequalities for $k\geq1$, it follows that
\begin{equation}
D=F(\mathbf{x}^1)\geq\delta\sum_{k=1}^\infty||\mathbf{x}^{k+1}-\mathbf{x}^k||_2^2,
\end{equation}
This in particular implies that $\lim\limits_{k\rightarrow\infty}(\mathbf{x}^{k+1}-\mathbf{x}^k)=\mathbf{0}$. Also $\{\mathbf{x}^k\}$ is bounded due to the condition (\textbf{C4}). Similar to the proof in the Theorem \ref{thm_con}, it is easy to show that any accumulation point of $\{\mathbf{x}^k\}$ is a stationary point $\mathbf{x}^*$ of problem (\ref{eq_mspro}). The convergence rate can be easily proved by summing (\ref{eq_proofPIREPS111}) for $k=1,\cdots,n$.
$\hfill\blacksquare$

\textbf{Remark:} The convergence results in the Theorem \ref{thm_ps_pro_ps} is similar to that in the Theorem \ref{thm_pro}. But the main difference of the proof is that we use the Pythagoras relation (\ref{eq_Pythagoras_ps}) instead of the Lipschitz continuous property (\ref{eq_lipp}) of the squared loss function $h(\mathbf{x}_1,\cdots,\mathbf{x}_S)=\frac{1}{2}\left\|\sum_{s=1}^S\mathbf{A}_s\mathbf{x}_s-\mathbf{b}\right\|_2^2$. The Pythagoras relation (\ref{eq_Pythagoras_ps}) is much tighter than (\ref{eq_lipp}).

%

\section{Convergence Analysis of PIRE-AU}

\begin{theo}\label{thm_ps_pro_ns}
For the loss function $h(\mathbf{x}_1,\cdots,\mathbf{x}_S)$ in problem (\ref{eq_mspro}), assume that $\nabla_sh(\mathbf{x}_1,\cdots,\mathbf{x}_S)$ is Lipschitz continuous with constant $L_s(h)$, $s=1,\cdots,S$. Let $D=F(\mathbf{x}^1)$, $\mu_s>\frac{L_s(h)}{2}$, $s=1,\cdots,S$, and $\delta=\min_s\{\mu_s-L_s(h)/2\}$. The sequence $\{\mathbf{x}^k\}$ generated by PIRE-AU in (\ref{eq_updatexns}) and (\ref{eq_updatewinpa}) satisfies the following properties:
\begin{enumerate}[(1)]
\item $F(\mathbf{x}^k)$ is monotonically decreasing, i.e. $F(\mathbf{x}^{k+1})\leq F(\mathbf{x}^k)$. Indeed,
\begin{equation*}
F(\mathbf{x}^k)-F(\mathbf{x}^{k+1})\geq\sum_{s=1}^S\left(\mu_s-\frac{L_s(h)}{2}\right)||\mathbf{x}_s^{k+1}-\mathbf{x}_s^k||^2\geq0;
\end{equation*}
\item The sequence $\{\mathbf{x}^k\}$ is bounded;
\item $\lim\limits_{k\rightarrow\infty}(\mathbf{x}^k-\mathbf{x}^{k+1})=\bm{0}$.
\item Any accumulation point of $\{\mathbf{x}^k\}$ is a stationary point $\mathbf{x}^*$ of problem (\ref{eq_mspro}).
\item For any $n\geq1$, $\min_{1\leq k\leq n}||\mathbf{x}^{k+1}-\mathbf{x}^{k}||_2^2\leq\frac{F(\mathbf{x}^1)-F(\mathbf{x}^*)}{n\delta}$.
\end{enumerate}
\end{theo}
\textbf{Proof}. Since $\mathbf{x}^{k+1}_s$ is the globally optimal solution to problem (\ref{eq_updatexns}), the zero vector is contained in the subgradient with respect to $\mathbf{x}_s$. That is, there exists $\mathbf{v}_s^{k+1}\in\partial\langle\mathbf{w}_s^k,\mathbf{g}(\mathbf{x}^{k+1}_s)\rangle$ such that
\begin{equation}\label{eq_proof1_ps}
\lambda\mathbf{v}^{k+1}_s+\nabla_sh(\mathbf{x}_1^{k+1},\cdots,\mathbf{x}_{s-1}^{k+1},\mathbf{x}_s^k,\cdots,\mathbf{x}_S^k)+\mu_s(\mathbf{x}^{k+1}_s-\mathbf{x}^{k}_s)=\mathbf{0}.
\end{equation}
A dot-product with $\mathbf{x}^{k+1}_s-\mathbf{x}^k_s$ on both sides of (\ref{eq_proof1_ps}) gives
\begin{equation}\label{eq_proof2_ps}
\begin{split}
\lambda\left\langle \mathbf{v}_s^{k+1},\mathbf{x}_s^{k+1}-\mathbf{x}_s^k\right\rangle+\left\langle \nabla_sh(\mathbf{x}_1^{k+1},\cdots,\mathbf{x}_{s-1}^{k+1},\mathbf{x}_s^k,\cdots,\mathbf{x}_S^k),\mathbf{x}_s^{k+1}-\mathbf{x}_s^k\right\rangle+\mu_s||\mathbf{x}_s^{k+1}-\mathbf{x}_s^k||_2^2=0.
\end{split}
\end{equation}
Recalling the definition of the subgradient of the convex function, we have
\begin{equation}\label{eq_nsproof1_ns}
\begin{split}
&\langle\mathbf{w}_s^k,\mathbf{g}(\mathbf{x}_s^k)-\mathbf{g}(\mathbf{x}_s^{k+1})\rangle\geq\langle\mathbf{v}_s^{k+1},\mathbf{x}_s^{k}-\mathbf{x}_s^{k+1}\rangle\\\
\end{split}
\end{equation}
Combining (\ref{eq_proof2_ps}) and (\ref{eq_nsproof1_ns}) gives
\begin{equation}\label{eq_nsproof1_ns22}
\begin{split}
\lambda\sum_{s=1}^S\langle\mathbf{w}_s^k,\mathbf{g}(\mathbf{x}_s^k)-\mathbf{g}(\mathbf{x}_s^{k+1})\rangle\geq\sum_{s=1}^S\left\langle\nabla_sh\left(\mathbf{x}_1^{k+1},\cdots,\mathbf{x}_{s-1}^{k+1},\mathbf{x}_s^k,\cdots,\mathbf{x}_S^k\right),\mathbf{x}_s^{k+1}-\mathbf{x}_s^{k}\right\rangle+\sum_{s=1}^S\mu_s||\mathbf{x}_s^{k+1} -\mathbf{x}_s^{k} ||_2^2.
\end{split}
\end{equation}
By the assumption (\textbf{C1}) that $f_s$ is concave, we get
\begin{equation}\label{eq_proof5_ns}
f_s(\mathbf{g}_s(\mathbf{x}_s^k))-f_s(\mathbf{g}_s(\mathbf{x}^{k+1}))\geq\langle\mathbf{w}_s^k,\mathbf{g}_s(\mathbf{x}_s^k)-\mathbf{g}_s(\mathbf{x}_s^{k+1})\rangle.
\end{equation}
Since $\nabla_sh(\mathbf{x}_1,\cdots,\mathbf{x}_S)$ is Lipschitz continuous, by using the property (\ref{eq_lipp}) for $s=1,\cdots,S$, we have
\begin{equation}\label{eq_nsproof2_ns}
\begin{split}
&h(\mathbf{x}_1^k,\cdots,\mathbf{x}_S^k)\\
\geq&h(\mathbf{x}_1^{k+1},\mathbf{x}_2^{k},\cdots,\mathbf{x}_S^k)-\left\langle\nabla_1h(\mathbf{x}_1^k,\cdots,\mathbf{x}_S^k),\mathbf{x}_1^{k+1}-\mathbf{x}_1^{k}\right\rangle-\frac{L_1(h)}{2}||\mathbf{x}_1^{k+1}-\mathbf{x}_1^{k}||_2^2\\
\geq&h(\mathbf{x}_1^{k+1},\mathbf{x}_2^{k+1},\mathbf{x}_3^k,\cdots,\mathbf{x}_S^k)-\left\langle\nabla_2h(\mathbf{x}_1^{k+1},\mathbf{x}_2^k,\cdots,\mathbf{x}_S^k),\mathbf{x}_2^{k+1}-\mathbf{x}_2^{k}\right\rangle-\left\langle\nabla_1h(\mathbf{x}_1^k,\cdots,\mathbf{x}_S^k),\mathbf{x}_1^{k+1}-\mathbf{x}_1^{k}\right\rangle\\
&-\frac{L_2(h)}{2}||\mathbf{x}_2^{k+1}-\mathbf{x}_2^{k}||_2^2-\frac{L_1(h)}{2}||\mathbf{x}_1^{k+1}-\mathbf{x}_1^{k}||_2^2\\
\geq&h(\mathbf{x}_1^{k+1},\cdots,\mathbf{x}_S^{k+1})-\sum_{s=1}^S\left\langle\nabla_sh(\mathbf{x}_{1}^{k+1},\cdots,\mathbf{x}_{s-1}^{k+1},\mathbf{x}_{s}^{k},\mathbf{x}_{s+1}^{k},\cdots,\mathbf{x}_{S}^{k}),\mathbf{x}_{s}^{k+1}-\mathbf{x}_{s}^{k}\right\rangle-\sum_{s=1}^S\frac{L_s(h)}{2}||\mathbf{x}_{s}^{k+1}-\mathbf{x}_{s}^{k}||_2^2.
\end{split}
\end{equation}
Combining (\ref{eq_nsproof1_ns22})(\ref{eq_proof5_ns}) and (\ref{eq_nsproof2_ns}) leads to
\begin{equation}
\begin{split}
&F(\mathbf{x}^k)-F(\mathbf{x}^{k+1})\\
=&\lambda\sum_{s=1}^S\left(f_s(\mathbf{g}_s(\mathbf{x}_s^k))-f_s(\mathbf{g}_s(\mathbf{x}_s^{k+1})) \right)+h(\mathbf{x}_1^k,\cdots,\mathbf{x}_S^k)-h(\mathbf{x}_1^{k+1},\cdots,\mathbf{x}_S^{k+1})\\
\geq&\sum_{s=1}^S\left(\mu_s-\frac{L_s(h)}{2}\right)||\mathbf{x}_{s}^{k+1}-\mathbf{x}_{s}^{k}||_2^2\\
\geq&\delta||\mathbf{x}^{k+1}-\mathbf{x}^{k}||_2^2\geq0.
\end{split}
\end{equation}
Thus $F(\mathbf{x}^k)$ is monotonically decreasing. Similarly, it is easy to see that $\lim\limits_{k\rightarrow\infty}(\mathbf{x}^k-\mathbf{x}^{k+1})=\bm{0}$. The sequence $\{\mathbf{x}^k\}$ is bounded due to condition (\textbf{C4}). This guarantees that $\{\mathbf{x}^k\}$ exists at least one accumulation point, and it is a stationary point as that in Theorem \ref{thm_con}. The convergence rate can be obtained in the same way.
 $\hfill\blacksquare$

\textbf{Remark:} The convergence results in the Theorem \ref{thm_ps_pro_ns} is similar to that in the Theorem \ref{thm_pro}. The main difference of the proof is that we use the Lipschitz continuous property for all each $\nabla_s h(\mathbf{x}_1,\cdots,\mathbf{x}_S)$, $s=1,\cdots,S$. This leads to (\ref{eq_nsproof2_ns}) which is tighter than (\ref{eq_proof6}), but for the multi-variable case.

\end{document}